\documentclass{article}

\usepackage[leqno]{amsmath}
\usepackage{amssymb}
\usepackage{proof}
\usepackage{fancybox}
\usepackage{enumerate}
\usepackage[dvips,ps,all]{xy}
\usepackage{fullpage}
\usepackage{theorem}
\usepackage{QED}
\usepackage{euscript}
\usepackage{stmaryrd}

\newenvironment{player}{\( \left\{\begin{array}{c}}{\end{array}\right\} \)}
\newenvironment{opponent}{\( \left(\begin{array}{c}}{\end{array}\right) \)}

\newcommand{\<}{\ensuremath{\langle}}
\renewcommand{\>}{\ensuremath{\rangle}}

\newcommand{\cat}[1]{\ensuremath{\mathbf{#1}}}

\newcommand{\script}[1]{\ensuremath{\EuScript {#1}}}

\newcommand{\bag}{\ensuremath{\mathrm{bag}}}

\newcommand{\x}{\ensuremath{\times}}

\newcommand{\ox}{\ensuremath{\varotimes}}
\newcommand{\bigox}{\ensuremath{\bigotimes}}
\newcommand{\ot}{\ensuremath{\varoplus}}
\newcommand{\bigot}{\ensuremath{\bigoplus}}
\newcommand{\vd}{\ensuremath{\vdash}}

\newcommand{\ra}{\ensuremath{\rightarrow}}

\newcommand{\Lra}{\ensuremath{\Longrightarrow}}

\newcommand{\pc}{\ensuremath{\Mapstochar=\!=\hspace{-1.2ex}=\Mapsfromchar}}
\newcommand{\ora}[1]{\ensuremath{\overrightarrow {#1}}}
\newcommand{\ola}[1]{\ensuremath{\overleftarrow {#1}}}

\newcommand{\hgt}{\ensuremath{\mathrm{hgt}}}
\newcommand{\cuthgt}{\ensuremath{\mathrm{cuthgt}}}

\newcommand{\bs}{\ensuremath{\backslash}}

\newcommand{\qqqquad}{\quad \quad \quad \quad}

\newcommand{\senta}[1]{\begin{center} #1 \end{center}}

\makeatletter

\newdimen\w@dth

\def\setw@dth#1#2{\setbox\z@\hbox{\scriptsize $#1$}\w@dth=\wd\z@
\setbox\@ne\hbox{\scriptsize $#2$}\ifnum\w@dth<\wd\@ne \w@dth=\wd\@ne \fi
\advance\w@dth by 1.2em}

\def\t@^#1_#2{\allowbreak\def\n@one{#1}\def\n@two{#2}\mathrel
{\setw@dth{#1}{#2}
\mathop{\hbox to \w@dth{\rightarrowfill}}\limits
\ifx\n@one\empty\else ^{\box\z@}\fi
\ifx\n@two\empty\else _{\box\@ne}\fi}}
\def\t@@^#1{\@ifnextchar_ {\t@^{#1}}{\t@^{#1}_{}}}

\def\t@left^#1_#2{\def\n@one{#1}\def\n@two{#2}\mathrel{\setw@dth{#1}{#2}
\mathop{\hbox to \w@dth{\leftarrowfill}}\limits
\ifx\n@one\empty\else ^{\box\z@}\fi
\ifx\n@two\empty\else _{\box\@ne}\fi}}
\def\t@@left^#1{\@ifnextchar_ {\t@left^{#1}}{\t@left^{#1}_{}}}

\def\two@^#1_#2{\def\n@one{#1}\def\n@two{#2}\mathrel{\setw@dth{#1}{#2}
\mathop{\vcenter{\hbox to \w@dth{\rightarrowfill}\kern-4ex
                 \hbox to \w@dth{\rightarrowfill}}%
       }\limits
\ifx\n@one\empty\else ^{\box\z@}\fi
\ifx\n@two\empty\else _{\box\@ne}\fi}}
\def\tw@@^#1{\@ifnextchar_ {\two@^{#1}}{\two@^{#1}_{}}}

\def\tofr@^#1_#2{\def\n@one{#1}\def\n@two{#2}\mathrel{\setw@dth{#1}{#2}
\mathop{\vcenter{\hbox to \w@dth{\rightarrowfill}\kern-4ex
                 \hbox to \w@dth{\leftarrowfill}}%
       }\limits
\ifx\n@one\empty\else ^{\box\z@}\fi
\ifx\n@two\empty\else _{\box\@ne}\fi}}
\def\t@fr@^#1{\@ifnextchar_ {\tofr@^{#1}}{\tofr@^{#1}_{}}}

\newdimen\W@dth
\def\setW@dth#1#2{\setbox\z@\hbox{$#1$}\W@dth=\wd\z@
\setbox\@ne\hbox{$#2$}\ifnum\W@dth<\wd\@ne \W@dth=\wd\@ne \fi
\advance\W@dth by 1.2em}

\def\T@^#1_#2{\allowbreak\def\N@one{#1}\def\N@two{#2}\mathrel
{\setW@dth{#1}{#2}
\mathop{\hbox to \W@dth{\rightarrowfill}}\limits
\ifx\N@one\empty\else ^{\box\z@}\fi
\ifx\N@two\empty\else _{\box\@ne}\fi}}
\def\T@@^#1{\@ifnextchar_ {\T@^{#1}}{\T@^{#1}_{}}}

\def\T@left^#1_#2{\def\N@one{#1}\def\N@two{#2}\mathrel{\setW@dth{#1}{#2}
\mathop{\hbox to \W@dth{\leftarrowfill}}\limits
\ifx\N@one\empty\else ^{\box\z@}\fi
\ifx\N@two\empty\else _{\box\@ne}\fi}}
\def\T@@left^#1{\@ifnextchar_ {\T@left^{#1}}{\T@left^{#1}_{}}}

\def\Tofr@^#1_#2{\def\N@one{#1}\def\N@two{#2}\mathrel{\setW@dth{#1}{#2}
\mathop{\vcenter{\hbox to \W@dth{\rightarrowfill}\kern-4ex
                 \hbox to \W@dth{\leftarrowfill}}%
       }\limits
\ifx\N@one\empty\else ^{\box\z@}\fi
\ifx\N@two\empty\else _{\box\@ne}\fi}}
\def\T@fr@^#1{\@ifnextchar_ {\Tofr@^{#1}}{\Tofr@^{#1}_{}}}

\def\Two@^#1_#2{\def\N@one{#1}\def\N@two{#2}\mathrel{\setW@dth{#1}{#2}
\mathop{\vcenter{\hbox to \W@dth{\rightarrowfill}\kern-4ex
                 \hbox to \W@dth{\rightarrowfill}}%
       }\limits
\ifx\N@one\empty\else ^{\box\z@}\fi
\ifx\N@two\empty\else _{\box\@ne}\fi}}
\def\Tw@@^#1{\@ifnextchar_ {\Two@^{#1}}{\Two@^{#1}_{}}}

\def\to{\@ifnextchar^ {\t@@}{\t@@^{}}}
\def\from{\@ifnextchar^ {\t@@left}{\t@@left^{}}}
\def\two{\@ifnextchar^ {\tw@@}{\tw@@^{}}}
\def\tofro{\@ifnextchar^ {\t@fr@}{\t@fr@^{}}}
\def\To{\@ifnextchar^ {\T@@}{\T@@^{}}}
\def\From{\@ifnextchar^ {\T@@left}{\T@@left^{}}}
\def\Two{\@ifnextchar^ {\Tw@@}{\Tw@@^{}}}
\def\Tofro{\@ifnextchar^ {\T@fr@}{\T@fr@^{}}}

\makeatother

\newtheorem{theorem}{Theorem}[section]
\newtheorem{lemma}[theorem]{Lemma}
\newtheorem{proposition}[theorem]{Proposition}
\newtheorem{corollary}[theorem]{Corollary}
{\theorembodyfont{\rmfamily} }
{\theorembodyfont{\rmfamily} \newtheorem{example}[theorem]{Example}}
{\theorembodyfont{\rmfamily} \newtheorem{remark}[theorem]{Remark}}

\newcommand{\channel}{\ensuremath{\mathsf{CProc}}}
\renewcommand{\phi}{\ensuremath{\varphi}}

\newcommand{\cone}{\ensuremath{\mathsf{c1}}}
\newcommand{\ctwo}{\ensuremath{\mathsf{c2}}}
\newcommand{\oned}{\ensuremath{\$1}}
\newcommand{\twod}{\ensuremath{\$2}}
\newcommand{\change}{\ensuremath{\mathsf{change}}}
\newcommand{\nochange}{\ensuremath{\mathsf{nochange}}}
\newcommand{\gal}{\ensuremath{\mathsf{gal}}}
\newcommand{\GAL}{\ensuremath{\mathsf{GAL}}}
\newcommand{\gum}{\ensuremath{\mathsf{gum}}}
\newcommand{\gumch}{\ensuremath{\mathsf{gumch}}}
\newcommand{\GUM}{\ensuremath{\mathsf{GUM}}}
\newcommand{\cc}{\ensuremath{\boldsymbol{:}}}

\title{\textbf{A language for multiplicative-additive linear logic}\thanks{
Department of Computer Science, University of Calgary, 2500 University
Drive NW, Calgary, Alberta, Canada T2N 1N4. Research partially supported by
NSERC, Canada. Diagrams were produced with the \Xy-pic package of K. Rose
and R. Moore and inferences with M. Tatsuya's \texttt{proof.sty}.}}
\author{J.R.B. Cockett and C.A. Pastro}
\date{\today}

\begin{document}
\maketitle

\begin{abstract}
A term calculus for the proofs in multiplicative-additive linear logic is 
introduced and motivated as a programming language for channel based 
concurrency. The term calculus is proved complete for a semantics in
linearly distributive categories with additives. It is also shown that proof
equivalence is decidable by showing that the cut elimination rewrites supply
a confluent rewriting system modulo equations.
\end{abstract}

\setcounter{section}{-1}
\section{Introduction} 

The purpose of this paper is to introduce a term logic for the
multiplicative-additive fragment of linear logic. Before introducing
this term logic, however, it is perhaps worth sharing with the reader how we
arrived at the language and why we think it is of some interest.
 
The idea of having a term logic for linear logic is certainly not new. In
fact, almost every researcher who has been heavily involved in trying to
understanding the proof theory of linear logic has found it necessary to
invent\footnote{Inevitably sometimes the word ``reinvent'' should be used
here. However, in defense of reinvention it should be said that, often, to
reinvent is the only way to really understand these languages.} such a
calculus. Perhaps the earliest attempt at a term logic for monoidal
categories was by Jay~\cite{jay89:languages} who essentially realized that
ordinary algebraic terms with no variable copying or elimination would do
the trick. More sophisticated attempts followed which linked the term
calculus to the proof theory of various fragments of linear logic, see,
e.g., Abramsky~\cite{abramsky93:computational} and Benton, Bierman, de
Paiva, and Hyland~\cite{benton93:termcal}.

An important component of Girard's approach to linear
logic~\cite{girard87:linear} was the introduction of proof nets. These, of
course, can also be regarded as a term logic in their own right and, indeed,
in~\cite{blute96:natural}, were explicitly introduced as such. There, after
straightening out Girard's one-sided proof nets into a two-sided form, they
were used as a basis for solving the coherence problem for the units. From
the point of view of a term logic, however, there is something peculiar
about these nets as their definition is not local: one has to
check a global correctness criterion before one can conclude that the net
is valid.  This condition arises as, in checking that the net corresponds
to a valid proof, one has to determine that there is a valid way of
assigning ``scopes'' to the inference rules. That this can be expressed as
a combinatoric condition on the nets was, of course, one of Girard's key
insights.

An interesting recent approach to providing a term logic was suggested by
Koh and Ong~\cite{koh98:type}. They realized that the tricky rewiring
conditions for the units which arose in \cite{blute96:natural} could be
expressed directly and quite clearly with scope changing rules. The first
author was very fortunate to have visited Koh and Ong in Oxford in 1996 and
to have had a chance to discuss this term logic with them. He was, of course,
particularly impressed by the fact that they had realized that this gave a
natural term logic for linearly distributive categories. It was clear that
they had a good idea.  However, their term logic never found any strong
resonance in the linear logic community. This was not really surprising: 
after all, the proof net technology and its correctness criterion had been
invented precisely to remove the necessity of keeping track of scopes. The
reintroduction of explicit scopes seemed like a step in the wrong direction
and made the utility of such an logic rather difficult to sell. To make
matters worse, the syntax of their term logic was concise to the point of
being cryptic: for an outsider the terms did not invite any particular
insight into their meaning.

While trying to sort out a process semantics for multiplicative-additive
linear logic (with both the additive and multiplicative units) we found
that it was very useful to have a term logic to express the processes.
The use of proof nets in linear logic makes it largely unnecessary to
have a term logic for the multiplicative fragment.  However, when one 
considers the multiplicatives together with the additives the value of a 
term logic becomes much more compelling. Even the minimalistic 
approach presented by Hughes and
van Glabbeek in \cite{hughes03:proof} cannot hide the fact that additive
proof nets are complicated combinatoric structures which are hard to
create in the way a programmer might create a program. 

In order to produce such a term logic we took a very literal translation of
the sequent calculus. When we showed this term logic to Robert Seely, and
explained the intended process semantics, he suggested ---while commenting
ironically on the importance of appearances rather than content in certain
academic circles--- that we should try presenting it in a programming
language style for processes. So we fiddled around with suitable keywords
and produced such a language and were horrified to find that we actually
liked the result! In particular, it had an immediate resonance with process
programming which made the terms almost readable to programmers.   

The point is that in presenting this term logic we are happy to claim almost
no originality. In fact, we would suggest that the strength of the language
is exactly that it has been plagiarized from a number of sources while 
slavishly following the sequent calculus. The scoping techniques should 
be quite recognizable from the term logic of Koh and Ong. Admittedly, we 
have unashamedly rearranged their terms to promote, following Robert's
excellent suggestion, a process reading. Finally, center stage is the syntax
for the additives which we were happy to simply borrow from the
$\pi$-calculus. The result is a term calculus which has the feel of a
programming language for processes.

Indeed, if you had been trying to put together a process language for
channel based communication you might well have come up with exactly the
same language. That the language is, in fact, the internal language for
linearly distributive categories with additives should be, perhaps, of more
than a passing interest.  

\subsection{Introduction to the term calculus}

As far as we know there has not actually been a proposal for a term calculus
for the multiplicative-additive fragment of linear logic. We should,
therefore, mention that although we borrow techniques very happily from
prior work there are some technical aspects which remain. For example,
the proof that term normalization modulo equations works is, for this
fragment, technically more challenging than in the purely multiplicative
fragment (compare our techniques to those in \cite{koh98:type}). While it is
not possible to include all the technical details in this paper the reader
should be aware that many of these details are available from the second
author's master's thesis~\cite{pastro:msc}. The subject of this thesis was
the complete additive fragment of linear logic\footnote{The complete
additive fragment includes, in particular, the additive units which not
surprisingly are more technically challenging to capture.} for which a full
and faithful process semantics was provided. This paper is part of the
continuation of that work, whose aim is to provide a full and faithful
process semantics for the complete multiplicative and additive fragments of
linear logic\footnote{The complete multiplicative and additive fragment
includes both the multiplicative units and the additive units. In particular, 
it should be mentioned we are fundamentally not assuming ``mix'', thus the
complications of handling the units in their full glory is present.}.

A slight warning to linear logicians: we do not assume that we have a
negation in the logic. This may seem to be an enormous omission, however,
those intimately familiar with the coherence issues of these settings will
know that, in fact, negation is a distraction whose presence or absence is
orthogonal to the real issues. If this sounds like a controversial statement
to the reader then we should perhaps also quietly mention that the initial
model does actually have negation (i.e., it is a $*$-autonomous category)
even though no negation is mentioned: in other words, this aspect, for the
initial model, comes along for the ride anyway. Thus, we are just being
scrupulously abstract in our approach and are therefore working at the 
level of linearly distributive categories with additives.

The objective of this paper is to highlight the term calculus we are using
and to thereby persuade you that it is quite reasonable to give a process
reading to the proofs of this fragment of linear logic. To this end we will
start with an old (and perhaps now unfashionable) example of Girard's
which involves obtaining a packet of Galois from a vending machine. We shall
show how to program it in our term logic.

We wish to describe the behavior of a vending machine which allows one to
select either smokers chewing gum or, for those that cannot quite kick the
habit, a packet of Galois cigarettes. The machine allows you to pay by
inserting a one or two dollar coin. A packet of Galois costs two dollars
while the gum only costs one dollar. There are four possible outcomes: 

\begin{enumerate}
\item You have inserted a dollar and have requested gum: you get a pack of gum.
\item You have inserted a dollar and requested a packet of Galois: your
dollar is returned.
\item You have inserted two dollars and requested some gum: you get a pack
of gum and a dollar in change.
\item You have inserted two dollars and requested a packet of Galois: you
get a packet of Galois.
\end{enumerate}
This means we wish to produce a process which is the proof of the following
sequent:
\[\alpha \cc \{\cone:\oned, \ctwo:\twod\} \ox \{\gal:\top,\gum:\top\} \vd
\beta \cc \{a:\GAL,b:\GUM \ox \{\change:\oned,\nochange:\top\},c:\oned\}
\]
The process starts with a single input channel $\alpha$ and a single output
channel $\beta$. Each channel has an attached protocol for the interaction
down that channel. The curly brackets denote coproduct types with components
which are tagged: the names used for the tags are called ``events'' and are 
sent and received as messages along the channels. The tensor operation 
$- \ox -$ allows one to bundle channels together.  

To model the functioning of the vending machine we need three functions
which transmute our resources:

\[\gal:\twod \ra \GAL, \quad \gum: \oned \ra \GUM, \quad \text{and} \quad
\gumch: \twod \ra \GUM \ox \oned.
\]

Thus, $\gal$ turns two dollars into a packet of Galois, $\gum$ turns one
dollar into a stick of gum, and $\gumch$ turns two dollars into a stick of
gum and one dollar of change. The fairness of the machine is guaranteed by
the fact that only these exchanges are allowed: to be correct we also have
to ensure (for example) that when the user inserts \$2 and asks for 
cigarettes he will not be given gum and a dollar change.

This process may be written as follows:

\texttt{\begin{tabbing}
sp\=lit $\alpha$ as $\alpha_1,\alpha_2$ in \\
\> in\=put $\alpha_1$ of \\
\>\> $\mid \cone$ \= $\mapsto$ input $\alpha_2$ of \\
\>\>\> $\mid \gal$ \= $\mapsto$ close $\alpha_2$ in \\
\>\>\>\> output $c$ on $\beta$ in $\alpha_1 = \beta$ \quad \% return coin \\
\>\>\> $\mid \gum$ $\mapsto$ close $\alpha_2$ in \\
\>\>\>\> ou\=tput $b$ on $\beta$ in \\
\>\>\>\>\> fo\=rk $\beta$ as \\
\>\>\>\>\>\> $\mid \beta_1$ with $\alpha_1 \mapsto \gum(\alpha_1;\beta_1)$
\quad \% return gum \\
\>\>\>\>\>\> $\mid \beta_2$ with $\quad \mapsto$ output $\nochange$ on
$\beta_2$ in end $\beta_2$ \\
\>\> $\mid \ctwo$ $\mapsto$ input $\alpha_2$ of \\
\>\>\> $\mid \gal$ $\mapsto$ close $\alpha_2$ in \\
\>\>\>\> output $a$ on $\beta$ in $\gal(\alpha_1;\beta)$ \quad \% return cig \\
\>\>\> $\mid \gum$ $\mapsto$ close $\alpha_2$ in \\
\>\>\>\> output $b$ on $\beta$ in \\
\>\>\>\>\> on $\gamma$ plug \\
\>\>\>\>\>\>  $\gumch(\alpha_2;\gamma)$ \quad \% return gum with change \\
\>\>\>\>\> to \\
\>\>\>\>\>\> sp\=lit $\gamma$ as $\gamma_1,\gamma_2$ in \\
\>\>\>\>\>\>\> fo\=rk $\beta$ as \\
\>\>\>\>\>\>\>\> $\mid \beta_1$ with $\gamma_1 \mapsto \gamma_1 = \beta_1$ \\
\>\>\>\>\>\>\>\> $\mid \beta_2$ with $\gamma_2 \mapsto$ output $\change$ on
$\beta_2$ in $\gamma_2 = \beta_2$
\end{tabbing}}

Input channels which are tensored can be used by a single process as they
come from independent sources: to communicate down the individual channels
one must \texttt{split} the channels. However, output channels when
tensored together have to be treated completely differently: they can be
dependent and so if they are used by the same process this may can cause
deadlock or livelock.  Thus, to use such channels one must \texttt{fork} the
current process into two independent processes.  When one forks a process in
this manner one must also decide which of the open communication channels
will be attached to which of the subprocesses.

One can only close a channel when the type is a unit ($\top$ or $\bot$).
When it is an input channel one can simply \texttt{close} the channel $\top$
and when it is an output one can \texttt{end} the channel provided there are
no other open channels.

The process above starts by splitting its one input channel, $\alpha$, 
into two distinct
channels. It then looks for an input event on the first channel: this is
provided by the user inserting either \$1 or \$2 into the machine. Next
the machine looks for input on the buttons. If the user has input \$1
and asked for a packet of Galois then we return his coin by selecting case
$(c)$ and connecting the input channel with the coin resource directly to an
output (so the coin is returned).

Let us consider what is perhaps the most complex case in the example:
namely, when \$2 is entered and the gum button is pressed. In this case we
start by choosing case $(b)$, then we want to transmute \$2 into gum and
\$1 (which are sent along channel $\gamma$). Next we want to pass these
things to the person using the machine so we split the resources. Here we
have to \texttt{fork} as the user may use these channels in such a way as
to make them dependent (for example he may always look for his change before
he takes his gum). We then independently return his change and provide the
gum in two different processes. To provide the change we send an event
$\change$ to indicate that change is due and then pass the change down that
channel.

Finally let us consider the case where \$1 is provided and $\gum$ is
selected. In this case we again have to fork the output into two processes
one of which returns the gum and the other of which handles the change. On
the latter we indicate that no change is due by sending the event
$\nochange$: this gives a tensor unit channel on output which
(provided that no channels are open) allows one to end the process.

\paragraph{Outline of this paper}

The main formal result of the paper is that there is a decision procedure
for this term logic. However, the fact that such a procedure exists is, we
feel, a fairly standard observation and not as important as establishing the 
term logic itself. To establish
the term logic we introduce it formally in Section~\ref{sec-term}. In fact,
we introduce it in two different syntactic flavors: the programming syntax,
as above, and a more succinct (and cryptic) representation to facilitate
the technical arguments.

In Section~\ref{sec-term-rewrites} we introduce a rewriting system on terms.
The rewriting corresponds to the cut elimination procedure and proof
equivalences for multiplicative-additive linear logic.

In Section~\ref{sec-polycat-sem} we show that the terms can be interpreted
in any linearly distributive category with sums and products and,
furthermore, form themselves into the free representable polycategory with 
additives. By looking at the maps in this polycategory we obtain the corresponding 
free linearly distributive category with additives.  These results are an extension of the
results in~\cite{pastro:msc} and use similar techniques.

In Section~\ref{sec-dp} we discuss the proof of categorical cut elimination.
This involves showing, in proof theoretic terms, that the cut elimination
procedure terminates (which is well-known). However, it also involves
showing the more demanding fact these rewrites are confluent modulo the
equations introduced by the permuting conversions. Once one has categorical
cut elimination the equality of proofs can be determined by using the
equations alone: the fact that the system is decidable follows from the
subformula property which limits the number of proofs and, thus, the number
of terms in an equivalence class.  

This is a somewhat unsatisfactory state of affairs in which to leave things,
as it does not indicate the complexity of the decision procedure.
It is possible to organize the decision procedure to be much more efficient. 
Unfortunately, even in this reorganization the procedure has some steps to 
accomplish which correspond to the classical rewiring problem discussed in 
\cite{blute96:natural} for the multiplicatives. Thus, the complexity
of this determination is dominated by the complexity of determining equality
in the multiplicative fragment.  The first author has conjectured (pessimistically) 
that this is exponential and this still remains an open problem. 

\section{Term calculi for multiplicative-additive linear logic}
\label{sec-term}

In this section we introduce two term calculi for multiplicative-additive
linear logic. We first recall the sequent rules for the multiplicative and
additive fragments of linear logic\footnote{We are using the symbols $+$
and $\x$ to refer to the categorical sum and product respectively, and the
symbols $\ox$ and $\ot$ for the categorical tensor and cotensor
respectively. Note then that this conflicts with Girard's notation
in~\cite{girard87:linear}, but is consistent with~\cite{blute96:natural}}:

\begin{center}
\ovalbox{
\parbox{63ex}{\begin{center}
\medskip
$\infer[\text{(identity)}]{A \vd A}{}$ \medskip
\[\begin{array}{ccc}
\infer[\text{(cotuple)}]{\hspace{-1ex} \Gamma,\sum\limits_{i \in I} X_i \vd
\Delta \hspace{1ex}}{\{\Gamma,X_i \vd \Delta\}_{i \in I}} &&
\infer[\text{(tuple)}]{\hspace{-1ex} \Gamma \vd \prod\limits_{i \in I} Y_i,
\Delta\hspace{1ex}}{\{\Gamma \vd Y_i,\Delta\}_{i \in I}}
\medskip\\
\infer[\text{(projection)}]{\Gamma,\prod\limits_{i \in I} X_i \vd \Delta}
{\Gamma,X_k \vd \Delta} &\quad&
\infer[\text{(injection)}]{\Gamma \vd \sum\limits_{i \in I} Y_i,\Delta}
{\Gamma \vd Y_k,\Delta}
\vspace{-1ex}
\end{array}\]
\senta{where $k \in I, I \not= \emptyset$}
\[\begin{array}{ccc}
\infer[\text{(ltensor)}]{\hspace{-1ex} \Gamma,\bigox\limits_{i \in I} X_i
\vd \Delta \hspace{1ex}}{\Gamma,\{X_i\}_{i \in I} \vd \Delta} &&
\infer[\text{(rpar)}]{\hspace{-1ex} \Gamma \vd \bigot\limits_{i \in I} Y_i,
\Delta \hspace{1ex}}{\Gamma \vd \{Y_i\}_{i \in I},\Delta}
\medskip\\
\infer[\text{(lpar)}]{\hspace{-1ex} \Gamma,\bigot\limits_{i \in I} X_i \vd
\Delta \hspace{1ex}}{\{\Gamma_i,X_i \vd \Delta_i\}_{i \in I}} &\quad&
\infer[\text{(rtensor)}]{\hspace{-1ex} \Gamma \vd \bigox\limits_{i \in I}
Y_i, \Delta \hspace{1ex}}{\{\Gamma_i \vd Y_i,\Delta_i\}_{i \in I}}
\vspace{-1ex}
\end{array}\]
\senta{where $\Gamma=[\Gamma_i]_{i \in I}$ and $\Delta=[\Delta_i]_{i \in
I}$}
$\infer[\text{(cut)}]{\Gamma,\Gamma' \vd \Delta,\Delta'}{\Gamma \vd
\Delta,Z & Z,\Gamma' \vd \Delta'}$
\medskip
\end{center}}}
\end{center}

The additive fragment consists of the cotuple, tuple, projection, and
injection rules along with the identity and cut rules. The multiplicative
fragment consists of the ltensor, rtensor, rpar, and lpar rules along with
the identity and cut rules.

This presentation uses operations indexed by arbitrary (finite) sets. We
will often write simply $i$ (e.g., $\sum_i X_i$) or $j$ and use it to mean
$i \in I$ or $j \in J$. Except in the injection and projection rules these
index sets are allowed to be empty: this gives the nullary rules. These
nullary rules are usually presented separately as:

\begin{center}
\ovalbox{\parbox{50ex}{
\medskip
\[\begin{array}{ccc}
\infer[\text{(cotuple)}]{\Gamma,\mathbf{0} \vd \Delta}{} &\quad&
\infer[\text{(tuple)}]{\Gamma \vd \mathbf{1},\Delta}{}
\medskip\\
\infer[\text{(ltensor)}]{\Gamma,\top \vd \Delta}{\Gamma \vd \Delta}
&\quad&
\infer[\text{(rpar)}]{\Gamma \vd \bot, \Delta}{\Gamma \vd \Delta}
\medskip\\
\infer[\text{(lpar)}]{\bot \vd}{} && \infer[\text{(rtensor)}] {\vd \top}{}
\end{array}
\]
}}\end{center}
We shall consider various augmentations of this basic logic:
 
\begin{itemize}
\item The ``initial logic'' is the logic with no atoms. Notice this is still
a non-trivial logic because of the symbols $\mathbf{0},\ \mathbf{1},\ \top$,
and $\bot$. We shall denote this logic as $\channel(\emptyset)$.
 
\item The ``pure logic'' is the logic with an arbitrary set of atoms A:
we shall denote this logic as $\channel(A)$.
 
\item The ``free logic'' is the logic with an arbitrary set of atoms and an
arbitrary set of non-logical axioms relating lists of atoms. For example, if
$f$ is a non-logical axiom relating $A,B$ to $C,D$, this may be denoted as
a morphism or an inference, i.e., as
\[f:A,B \ra C,D \qquad \text{or} \qquad \infer{f :: A,B \vd C,D}{}
\]
The atoms will be regarded as objects in a polycategory and the axioms as
maps in that polycategory (with the ``essential cuts'' being provided by
composition). If the polycategory is denoted \cat{A}, the resulting logic
will be denoted $\channel(\cat{A})$.
\end{itemize}
 
If we think of the atoms of a pure logic as forming a discrete category (a
category where the only morphisms are the identities), the free logic on
this discrete category is just the ``pure'' logic. Each variant above
therefore includes the previous variant, and as it is more general, we
shall tend to only consider the free logic.  It is worth noting that 
a simple inductive argument shows that the logic will have negation whenever
the polycategory \cat{A} has negation.

\subsection{Formulas as protocols} \label{sec-annform}

The term representations are dependent on formulas annotated with ``channel
names'' and ``events'' which are derived from the process semantics view in 
which the formulas represent protocols. The reader interested in a more 
complete story is referred to~\cite{pastro:msc}.

Each formula is assigned a channel name, which we denote using Greek
letters, e.g., if $X$ is a formula, $\alpha \cc X$ is a formula attached 
to channel $\alpha$. A formula is just formula of linear logic (without
negation) which is presented using a tagged notation:

\begin{itemize}
\item if $X$ is atomic then it is left as is.
\item $X = \sum_i X_i$ is tagged with events as $\{a_i:X_i | i \in I\}$,
    where $a_i \neq a_j$ when $i \neq j$.
\item $X = \prod_i X_i$ is tagged with events as $(a_i:X_i| i \in I)$,
    where $a_i \neq a_j$ when $i \neq j$.
\item $X = \bigox_i X_i$ is tagged with channel names as $\bigox_i \alpha_i
    \cc X_i$, where $\alpha_i \neq \alpha_j$ when $i \neq j$.
\item $X = \bigot_i X_i$ is tagged as $\bigot_i \alpha_i \cc X_i$, where
    $\alpha_i \neq \alpha_j$ when $i \neq j$.
\end{itemize}
For example, the linear logic formula $\{W,X\} \ox (Y \ot Z)$ in this tagged
notation might become:
\[\alpha \cc \{a:W,b:X\} \ox \beta \cc (\beta_1 \cc Y \ot \beta_2 \cc Z)
\]

The channel names occurring as the tags for the formulas in a sequent must
be distinct as they are used as references. To simplify our term
representation and to avoid channel name conversions we will use distinct
channel names within the formulas as well. Furthermore, a convention we
adopt is to assign implicit tags to multiplicatives: thus, $\alpha \cc
\bigox_i X_i$ or $\alpha \cc \bigot_i X_i$ have implicitly assigned their
constituents the channel names $\alpha_i$, for $i \in I$. 

We are now ready to introduce the term calculi. Two different syntaxes will
be used to represent proofs in this system: a ``programming language''
syntax and a more succinct term representation. We begin with the
programming syntax.

\subsection{Programming language syntax} \label{sec-plrep}

The \textbf{programming syntax} provides an explicit process reading for
the proofs. Unfortunately, while this representation is quite intuitive, it
is also rather verbose which is why we will introduce a more succinct
representation in the next section.

The term formation rules for the programming language are given in
Table~\ref{tab-plrep}. The language is self dual and so only half the rules
are presented. To reduce the overload strain on colons we will use ``::''
to denote the term-type membership relation, e.g., $t::U \vd V$ will mean
that $t$ is a term of type $U \vd V$, where $U$ (say) may be of the form
$a:X$. Note also that we introduce special syntax (\texttt{stop},
\texttt{close}, and \texttt{end}) for the nullary cases.

The notation $\Lambda(\Gamma \vd \Delta)$ represents the set of domain and 
codomain  channels of the sequent, i.e., the set of channels tags of the
formulas in the sequent. For example,

\[\Lambda(\alpha \cc \{a:W,b:X\}, \alpha' \cc A \vd \beta \cc (Y \ot Z)) 
= \{\alpha,\alpha',\beta\}
\]
Note that $\Lambda$ does not extract any internal channel names.

\begin{table}
\begin{center}
\ovalbox{
\parbox[c]{75ex}{
\[\begin{array}{c}
\infer{\alpha \equiv_A \beta :: \alpha \cc A \vd \beta \cc A}{}
\bigskip\\
\infer{\texttt{input on } \alpha \texttt{ of } \mid_i a_i \mapsto f_i ::
\Gamma \vd \alpha \cc \prod_{i} a_i:X_i,\Delta}
{\{f_i::\Gamma \vd \alpha \cc X_i, \Delta\}_{i}}
\bigskip\\
\infer{\texttt{output } a_k \texttt{ on } \alpha \texttt{ in } f ::
\Gamma \vd \alpha \cc \sum_{i} a_i:X_i,\Delta}{f::\Gamma \vd \alpha \cc
X_k,\Delta} \smallskip\\
\text{where } k \in I, I \neq \emptyset
\bigskip\\
\infer{\texttt{split } \alpha \texttt{ as } (\alpha_i)_{i}
\texttt{ in } f:: \Gamma \vd \alpha \cc \bigot_i \alpha_i \cc X_i, \Delta}
{f::\Gamma \vd \{\alpha_i \cc X_i\}_i, \Delta}
\bigskip\\
\infer{\texttt{fork } \alpha \texttt{ as } \mid_i \alpha_i \texttt{ with }
\Lambda(\Gamma_i,\Delta_i) \mapsto f_i :: \Gamma \vd \alpha \cc \bigox_{i}
\alpha_i \cc X_i,\Delta}{\{f_i :: \Gamma_i \vd \alpha_i \cc X_i,\Delta_i\}_i}
\smallskip\\
\text{where } \Gamma = [\Gamma_i]_i \text{ and } \Delta = [\Delta_i]_i
\bigskip\\
\infer{\texttt{on } \gamma \texttt{ plug } f \texttt{ to } g ::
\Gamma,\Gamma' \vd \Delta,\Delta'}{f::\Gamma \vd \Delta,\gamma \cc Z &
g:: \gamma \cc Z,\Gamma' \vd \Delta'}
\end{array}\]
Notation for nullary cases:
\[\begin{array}{c}
\infer{\texttt{stop } \alpha :: \Gamma \vd \alpha \cc \mathbf{1},\Delta}{}
\bigskip\\
\infer{\texttt{close } \alpha \texttt{ in } f :: \Gamma \vd \alpha \cc
\bot,\Delta}{f::\Gamma \vd \Delta}
\bigskip\\
\infer{\texttt{end } \alpha ::\ \vd \alpha \cc \top}{}
\end{array}\]}}
\end{center}
\caption{``Programming language'' term formation rules}
\label{tab-plrep}
\end{table}

A non-logical axiom introduced as $f:\alpha_1 \cc X_1,\ldots, \alpha_m \cc
X_m \ra \beta_1 \cc Y_1,\ldots, \beta_n \cc Y_n$ will be represented by the
term $f(\alpha_1,\ldots,\alpha_m; \beta_1,\ldots,\beta_n)$. The identity
map on atoms, $\alpha \cc X \ra \beta \cc X$, can then be represented
using this notation as $1_X(\alpha;\beta)$, although we will prefer the
notation $\alpha \equiv_X \beta$. Two terms are \textbf{channel equivalent}
if they are equivalent up to channel name conversion: clearly the renaming
of channels does not affect the meaning of the terms.

\begin{example} \label{exam-plr}
The following proof of the linear distribution map

\[\infer{\alpha \cc A \ox (B \ot C) \vd \beta \cc B \ot (A \ox C)}{
\infer{\alpha_1 \cc A, \alpha_2 \cc B \ot C \vd \beta_1 \cc B, \beta_2 \cc
A \ox C}{
\infer{\alpha_{21} \cc B \vd \beta_1 \cc B}{} &
\infer{\alpha_1 \cc A,\alpha_{22} \cc C \vd \beta_2 \cc A \ox C}{
\infer{\alpha_1 \cc A \vd \beta_{21} \cc A}{} & \infer{\alpha_{22} \cc C \vd
\beta_{22} \cc C}{}}}}
\]
is represented in the programming syntax as follows:
\texttt{\begin{tabbing}
spl\=it $\alpha$ as $\alpha_1,\alpha_2$ in \\
\> spl\=it $\beta$ as $\beta_1,\beta_2$ in \\
\>\> for\=k $\alpha_2$ as \\
\>\>\> $\mid \alpha_{21}$ with $\beta_1 \mapsto \alpha_{21} \equiv_B
\beta_1$ \\
\>\>\> $\mid \alpha_{22}$ with $\alpha_1,\beta_2 \mapsto$ for\=k $\beta_2$
as \\
\>\>\>\> $\mid \beta_{21}$ with $\alpha_1$ in $\alpha_1 \equiv_A \beta_{21}$
\\
\>\>\>\> $\mid \beta_{22}$ with $\alpha_{22}$ in $\alpha_{22} \equiv_C
\beta_{22}$
\end{tabbing}}
\end{example}

\subsection{A term calculus representation} \label{sec-tcrep}

In this section we introduce a more succinct representation for proofs
which we call the \textbf{term calculus}. The term formation rules are
given in Table~\ref{tab-tcr}. Here we provide all the rules, however, again
notice that this representation is also self dual.

\begin{table}
\begin{center}
\ovalbox{
\parbox[c]{105ex}{\begin{center} \medskip
$\infer{\alpha \equiv_A \beta :: \alpha \cc A \vd \beta \cc A}{}$
\[\begin{array}{ccc}
\infer{\alpha\{a_i \mapsto f_i\}_i :: \Gamma, \alpha \cc \sum_i a_i:X_i
\vd \Delta}{\{f_i :: \Gamma,\alpha \cc X_i \vd \Delta\}_i}
&\qquad&
\infer{\alpha\{a_i \mapsto f_i\}_{i} :: \Gamma \vd \alpha \cc \prod_i
a_i:X_i,\Delta}{\{f_i :: \Gamma \vd \alpha \cc X_i, \Delta\}_i}
\bigskip\\
\infer{\alpha[a_k] f :: \Gamma, \alpha \cc \prod_i a_i:X_i \vd \Delta}
{f :: \Gamma,\alpha \cc X_k \vd \Delta} &&
\infer{\alpha[a_k] f :: \Gamma \vd \alpha \cc \sum_i a_i:X_i,\Delta}
{f :: \Gamma \vd \alpha \cc X_k,\Delta}
\vspace{-1ex}
\end{array}\]
\senta{where $k \in I, I \neq \emptyset$}
\[\begin{array}{ccc}
\infer{\alpha\<(\alpha_i)_i \mapsto f\> :: \Gamma, \alpha \cc \bigox_i
\alpha_i \cc X_i \vd \Delta}{f :: \Gamma,\{\alpha_i \cc X_i\}_i \vd \Delta} &&
\infer{\alpha\<(\alpha_i)_i \mapsto f\>:: \Gamma \vd \alpha \cc \bigot_i
\alpha_i \cc X_i, \Delta}{f::\Gamma \vd \{\alpha_i \cc X_i\}_i,\Delta}
\bigskip\\
\infer{\alpha\<\alpha_i \mid \Lambda(\Gamma_i,\Delta_i) \mapsto f_i\>_i ::
\Gamma,\alpha \cc \bigot_i \alpha_i \cc X_i \vd \Delta}
{\{f_i :: \Gamma_i,\alpha_i \cc X_i \vd \Delta_i\}_i} &&
\infer{\alpha\<\alpha_i \mid \Lambda(\Gamma_i,\Delta_i) \mapsto f_i\>_i ::
\Gamma \vd \alpha \cc \bigox_i \alpha_i \cc X_i,\Delta}
{\{f_i :: \Gamma_i \vd \alpha_i \cc X_i,\Delta_i\}_i}
\vspace{-1ex}
\end{array}\]
\senta{where $\Gamma = [\Gamma_i]_i$ and $\Delta = [\Delta_i]_i$}
$\infer{f ;_\gamma g :: \Gamma,\Gamma' \vd \Delta,\Delta'}{f :: \Gamma \vd
\Delta,\gamma \cc Z & g :: \gamma \cc Z,\Gamma' \vd \Delta'}$
\medskip
\end{center}}}
\end{center}
\caption{Term calculus formation rules}
\label{tab-tcr}
\end{table}

In this representation (as well as in the programming language
representation) the notation does not directly indicate on which side of
the turnstile an active channel sits. This can be inferred from the term
and also by inspecting the type. However, it is sometimes useful to make
this information more explicit and to do this we will indicate domain
channels by an overarrow pointing left, $\ola{\alpha}$, and codomain
channels by an overarrow pointing right, $\ora{\alpha}$.

\begin{example}
The following is the term calculus representation of the proof found in
Example~\ref{exam-plr} above:

\[\ola{\alpha}\<(\alpha_1,\alpha_2) \mapsto
\ora{\beta}\<(\beta_1,\beta_2) \mapsto \ola{\alpha}_2\left\<
\begin{array}{ll}
\alpha_{21} \mid \ora{\beta}_1 & \mapsto \ola{\alpha}_{21} \equiv_B \ora{\beta}_1 \smallskip\\
\alpha_{22} \mid \ola{\alpha}_1,\ora{\beta}_2 & \mapsto \ora{\beta}_2
  \left\<\begin{array}{ll}
  \beta_{21} \mid \ola{\alpha}_1 & \mapsto \ola{\alpha}_1 \equiv_A \ora{\beta}_{21} \\
  \beta_{22} \mid \ola{\alpha}_{22}  & \mapsto \ola{\alpha}_{22} \equiv_C \ora{\beta}_{22}
  \end{array}\right\>
\end{array}\right\>\>\>
\]
\end{example}

It is clear that the term calculus is more succinct and will be more convenient  
for manipulating the terms, accordingly we shall use this representation
for the remainder of the paper.

\section{Term rewrites and equivalences} \label{sec-term-rewrites}

Two terms $f$ and $g$ may be composed on a channel $\gamma$ when $\gamma$
is a codomain channel of $f$ of type $X$ and a domain channel of $g$ of the
same type.  Furthermore, $\gamma$ must be the only channel name $f$ and $g$ 
have in common.  That is, the two terms must be of the form $f::\Gamma_1
\ra \Gamma_2,\gamma \cc X$ and $g::\gamma \cc X, \Delta_1 \ra \Delta_2$,
where $\gamma$ is the only channel name in common. Composing these terms on
$\gamma$ is denoted:
\[f;_\gamma g :: \Gamma_1,\Delta_1 \ra \Gamma_2,\Delta_2
\]
The requirement that $\gamma$ is the only channel name in common ensures
that after composition all the channel names will be distinct. In general 
this means that in order to perform a composition there may be need to
perform a channel name conversion to arrange that the two terms have this
form.

\subsection{Cut elimination rewrites} \label{sec-cut-elim}

Composition is exactly the cut rule. Thus, the dynamics of composition 
corresponds precisely to the cut elimination process. We describe the cut
elimination rewrites below. As our term calculus does not distinguish between domain 
and codomain channels dual proof rewrites written as terms will be identical. 
However, we shall be careful to indicate how the inference rules give rise to each 
of the rewrites.

The set of unbound channels in a term $t$ will be denoted $\Lambda(t)$ and is the 
same as $\Lambda(\Gamma \vd \Delta)$ where $t::\Gamma \vd \Delta$.

The rewrites are as follows:

\begin{itemize}
\item Identity-sequent (sequent-identity):

$\begin{array}{lrcl}
(1) & f ;_\gamma 1 &\Lra& f \medskip\\
(2) & 1 ;_\gamma f &\Lra& f 
\end{array}$

\item Cotuple-sequent (sequent-tuple), tuple-sequent (sequent-cotuple):

$\begin{array}{lrcl}
(3) & \alpha\{a_i \mapsto f_i\}_i ;_\gamma g &\Lra&
\alpha\{a_i \mapsto f_i ;_\gamma g\}_i \medskip\\
(4) & g ;_\gamma \alpha\{a_i \mapsto f_i\}_i &\Lra&
\alpha\{a_i \mapsto g ;_\gamma f_i\}_i
\end{array}$

\item Injection-sequent (sequent-projection), projection-sequent
(sequent-injection):

$\begin{array}{lrcl}
(5) & \alpha[a] f ;_\gamma g &\Lra& \alpha[a] (f ;_\gamma g) \medskip\\
(6) & g ;_\gamma \alpha[a] f &\Lra& \alpha[a] (g ;_\gamma f)
\end{array}$

\item Ltensor-sequent (sequent-rpar), rpar-sequent (sequent-ltensor):

$\begin{array}{lrcl}
(7) & \alpha\<(\alpha_i)_i \mapsto f\> ;_\gamma g &\Lra&
\alpha\<(\alpha_i) \mapsto f ;_\gamma g\> \medskip\\
(8) & g ;_\gamma \alpha\<(\alpha_i)_i \mapsto f\> &\Lra&
\alpha\<(\alpha_i) \mapsto g ;_\gamma f\>
\end{array}$

\item Lpar-sequent (sequent-rtensor), rtensor-sequent (sequent-lpar): here
we suppose $\gamma \in \Lambda_k$.

$\begin{array}{lrcl}
(9) & \alpha\<\alpha_i \mid \Lambda_i \mapsto f_i\>_i ;_\gamma g &\Lra&
\alpha\left\<\begin{array}{l}
\alpha_i \mid \Lambda_i \mapsto f_i \\
\alpha_k \mid (\Lambda_k \cup \Lambda(g)) \setminus \{\gamma\} \mapsto
f_k ;_\gamma g
\end{array}\right\>_{i \neq k} \medskip\\
(10) & g ;_\gamma \alpha\<\alpha_i \mid \Lambda_i \mapsto f_i\>_i &\Lra&
\alpha\left\<\begin{array}{l}
\alpha_i \mid \Lambda_i \mapsto f_i \\
\alpha_k \mid (\Lambda_k \cup \Lambda(g)) \setminus \{\gamma\} \mapsto
g ;_\gamma f_k 
\end{array}\right\>_{i \neq k}
\end{array}$

\item Injection-cotuple (tuple-projection):

$\begin{array}{lrcl}
(11) & \gamma[a_k] f ;_\gamma \gamma\{a_i \mapsto g_i\}_i &\Lra&
f ;_\gamma g_k \medskip\\
(12) & \gamma\{a_i \mapsto g_i\}_i ;_\gamma \gamma[a_k] f  &\Lra&
g_k ;_\gamma f
\end{array}$

\item Rtensor-ltensor (rpar-lpar):

$\begin{array}{lrcl}
(13) & \gamma\<\alpha_i \mid \Lambda_i \mapsto f_i\>_{i \in \{1,\ldots,n\}}
;_\gamma \gamma\<(\alpha_i)_{i \in \{1,\ldots,n\}} \mapsto g\> &\Lra&
f_n ;_{\alpha_n} ( \cdots (f_2 ;_{\alpha_2} (f_1 ;_{\alpha_1} g)) \cdots)
\medskip\\
(14) &\gamma\<(\alpha_i)_{i \in \{1,\ldots,n\}} \mapsto g\> ;_\gamma 
\gamma\<\alpha_i \mid \Lambda_i \mapsto f_i\>_{i \in \{1,\ldots,n\}}
&\Lra&
(( \cdots (g ;_{\alpha_n} f_n) \cdots );_{\alpha_2} f_2) ;_{\alpha_1} f_1
\end{array}$
\end{itemize}

In order to obtain a normal form for sequent derivations, and hence terms,
we would like to show that the cut elimination rewrites are Church-Rosser.
However, this is easily seen to not be the case. Consider a derivation with
a projection and an injection immediately above the cut:
\[\infer{\Gamma_1,\prod_i X_i,\Delta_1 \vd \Gamma_2,\sum_j Y_j,\Delta_2}{
\infer{\Gamma_1,\prod_i X_i \vd \Gamma_2,Z}{
\infer{\Gamma_1,X_k \vd \Gamma_2,Z}{\pi}} &
\infer{Z,\Delta_1 \vd \sum_j Y_j,\Delta_2}{
\infer{Z,\Delta_1 \vd Y_l,\Delta_2}{\pi'}}}
\]
In this case one may apply the projection-sequent rewrite or the
sequent-injection rewrite to reduce the derivation, but there seems to be no
way in which to resolve this pair. This motivates the use of additional
rewrites which will allow us to interchange these rules (and all the other
critical pairs). These rewrites (which we are denoting by $\pc$) are called
the \textbf{permuting conversions} and are as follows:

\begin{itemize}
\item Cotuple-cotuple (tuple-tuple), cotuple-tuple interchange:

$\begin{array}{lrcl}
(15) & \alpha\{a_i \mapsto \beta\{b_j \mapsto f_{ij}\}_j\}_i &\pc&
\beta\{b_j \mapsto \alpha\{a_i \mapsto f_{ij}\}_i\}_j
\end{array}$

\item Cotuple-injection (projection-tuple), cotuple-projection
(injection-tuple) interchange:

$\begin{array}{lrcl}
(16) & \alpha\{a_i \mapsto \beta[b] f_i\}_i &\pc&
\beta[b] \alpha\{a_i \mapsto f_i\}_i
\end{array}$

\item Cotuple-ltensor (rpar-tuple), cotuple-rpar (ltensor-tuple) interchange:

$\begin{array}{lrcl}
(17) & \alpha\{a_i \mapsto \beta\<(\beta_j)_j \mapsto f_i\>\}_i &\pc&
\beta\<(\beta_j)_j \mapsto \alpha\{a_i \mapsto f_i\}_i\>
\end{array}$

\item Cotuple-rtensor (lpar-tuple), cotuple-lpar (rtensor-tuple) interchange:
here we suppose $\alpha \in \Lambda_k$.

$\begin{array}{lrcl}
(18) & \alpha\{a_i \mapsto \beta\left\<\begin{array}{l}
\beta_j \mid \Lambda_j \mapsto g_j \\
\beta_k \mid \Lambda_k \mapsto f_i
\end{array}\right\>_{j \neq k}\}_i
&\pc&
\beta\left\<\begin{array}{l}
\beta_j \mid \Lambda_j \mapsto g_j \\
\beta_k \mid \Lambda_k \mapsto \alpha\{a_i \mapsto f_i\}_i
\end{array}\right\>_{j \neq k}
\end{array}$

\item Injection-injection (projection-projection), injection-projection
interchange:

$\begin{array}{lrcl}
(19) & \alpha[a] (\beta[b] f) &\pc& \beta[b] (\alpha[a] f)
\end{array}$

\item Injection-ltensor (rpar-projection), injection-rpar
(ltensor-projection) interchange:

$\begin{array}{lrcl}
(20) & \alpha[a] \beta\<(\beta_j)_j \mapsto f\> &\pc&
\beta\<(\beta_j)_j \mapsto \alpha[a] f\>
\end{array}$

\item Injection-rtensor (lpar-projection), injection-lpar (rtensor-projection)
interchange: here we suppose $\alpha \in \Lambda_k$.

$\begin{array}{lrcl}
(21) & \alpha[a] \beta\left\<\begin{array}{l}
\beta_j \mid \Lambda_j \mapsto g_j \\
\beta_k \mid \Lambda_k \mapsto f
\end{array}\right\>_{j \neq k}
&\pc&
\beta\left\<\begin{array}{l}
\beta_j \mid \Lambda_j \mapsto g_j \\
\beta_k \mid \Lambda_k \mapsto \alpha[a] f
\end{array}\right\>_{j \neq k}
\end{array}$

\item Ltensor-ltensor (rpar-rpar), ltensor-rpar interchange:

$\begin{array}{lrcl}
(22) & \alpha\<(\alpha_i)_i \mapsto \beta\<(\beta_j)_j \mapsto f\>\> &\pc&
\beta\<(\beta_j)_j \mapsto \alpha\<(\alpha_i)_i \mapsto f\>\>
\end{array}$

\item Ltensor-rtensor (lpar-rpar), ltensor-lpar (rtensor-rpar) interchange:

$\begin{array}{lrcl}
(23) & \alpha\left\<(\alpha_i)_i \mapsto \beta\left\<\begin{array}{l}
\beta_j | \Lambda_j \mapsto g_j \\
\beta_k | \Lambda_k \cup (\bigcup_i \{\alpha_i\}) \mapsto f
\end{array}\right\>_{j \neq k}\right\>
\hspace{-1ex} &\pc& \hspace{-0.5ex}
\beta\left\<\begin{array}{l}
\beta_j | \Lambda_j \mapsto g_j \\
\beta_k | \Lambda_k \cup \{\alpha\} \mapsto \alpha\<(\alpha_i)_i \mapsto f\>
\end{array}\right\>_{j \neq k}
\end{array}$

\item Rtensor-rtensor (lpar-lpar), rtensor-lpar interchange: let
$\mathbf{\Lambda} = \{\beta\} \cup \Lambda_k \cup (\bigcup_j \Lambda'_j)$ and
$\mathbf{\Lambda'} = \{\alpha\} \cup \Lambda'_{k'} \cup (\bigcup_i \Lambda_i)$.

$\begin{array}{lrcl}
(24) & \alpha\left\<\begin{array}{l}
\alpha_i | \Lambda_i \mapsto f_i \\
\alpha_k | \mathbf{\Lambda} \mapsto
  \beta\left\<\begin{array}{l}
  \beta_j | \Lambda'_j \mapsto g_j \\
  \beta_{k'} | \Lambda'_{k'} \mapsto h
  \end{array}\right\>_{j \neq k'}
\end{array}\right\>_{i \neq k}
\hspace{-4ex} &\pc& \hspace{-0.5ex}
\beta \left\<\begin{array}{l}
\beta_j | \Lambda'_j \mapsto g_j \\
\beta_{k'} | \mathbf{\Lambda'} \mapsto
  \alpha\left\<\begin{array}{l}
  \alpha_i | \Lambda_i \mapsto f_i \\
  \alpha_k | \Lambda_k \mapsto h
  \end{array}\right\>_{i \neq k}
\end{array}\right\>_{j \neq k'}
\end{array}$
\end{itemize}

\section{Polycategorical semantics} \label{sec-polycat-sem}

The reduction rules and the permuting conversions together define an
equivalence relation (which we denote by $\sim$) on the derivations of a
sequent. This allows us to form for each sequent a hom-set consisting of 
terms modulo this equivalence (based on a starting polycategory \cat{A}).  
The cut operation then provides polycategorical composition.  

This delivers a polycategory which we shall denote \channel(\cat{A}). The
purpose of this section is to prove that this polycategory has additives 
and is representable.  This amounts to saying that its maps form a linearly
distributive category with additives (see \cite{cockett99:linearly}).

\begin{theorem} \label{prop-cha-poly}
\channel(\cat{A}) is a representable additive polycategory whose objects are 
the formulas of the logic, and whose morphisms are $\sim$-equivalence classes 
of derivations.
\end{theorem}

The fact that the terms modulo equivalence provide a polycategory is proved
in Appendix~\ref{sec-proof-polycat}. The proof consists of providing
identity morphisms and showing that cut is associative and satisfies the
interchange law. Below we discuss the remaining structure:

\subsection{Additives} \label{sec-sums}

Our next task is to establish that the polycategory $\channel(\cat{A})$
has additives, i.e., poly-sums and poly-products. In a polycategory \cat{A}
an object $\sum_i X_i \in \cat{A}$ is said to be the \textbf{poly-sum} (or
\textbf{poly-coproduct}) of a family of objects $X_i \in \cat{A}$, for $i
\in I$, in case there is a natural bijective correspondence

\[\infer={\alpha\{f_i\}_i ::\Gamma,\alpha \cc \sum_i X_i \ra \Delta}
{\{f_i::\Gamma,\alpha \cc X_i \ra \Delta\}_i}
\]
where, by natural, it is meant that the equivalences
\[h ;_\gamma \alpha\{f_i\}_i = \alpha\{h ;_\gamma f_i\}_i
\qquad \text{and} \qquad
\alpha\{f_i\}_i ;_\gamma h  = \alpha\{f_i ;_\gamma h\}_i
\]
hold (when $\alpha \neq \gamma$). Poly-products in polycategories are
(as might be expected) dual to poly-coproducts.

\begin{proposition} \label{prop-hassums}
$\channel(\cat{A})$ has additives.
\end{proposition}

\begin{proof}
In order to establish that $\channel(\cat{A})$ has finite sums and finite
products it must be shown that the correspondences

\[\vcenter{\infer={\Gamma,\alpha \cc \sum_i i:X_i \ra \Delta}
{\{\Gamma,\alpha \cc X_i \ra \Delta\}_i}}
\qquad \text{and} \qquad
\vcenter{\infer={\Gamma \ra \beta \cc \prod_j j:Y_j,\Delta}
{\{\Gamma \ra \beta \cc Y_j,\Delta\}_j}}
\]
are bijective and natural. We will prove the case for sums, products are
handled dually.

Suppose we have maps $t::\Gamma,\alpha \cc \sum_i i:X_i \ra \Delta$ and
$s_i::\Gamma,\alpha \cc X_i \ra \Delta$, for $i \in I$. The collection of
maps $\{s_i\}_i$ can, via the cotupling derivation, be used to construct a map

\senta{$\psi(\{s_i\}_i) =
\alpha\{i \mapsto s_i\}_i :: \Gamma,\alpha \cc \sum_i i:X_i \ra \Delta$}
and $t$ may be cut with $\alpha[k]1_{X_k} :: X_i \ra \alpha \cc
\sum_i X_i$, for $k \in I$, to get a map

\[\phi_k(t) = \alpha[k]1_{X_k} ;_\alpha t :: \Gamma,X_k \ra \Delta
\]

To prove that this correspondence is bijective we must establish both
$\phi_k(\psi(\{s_i\}_i)) = s_k$, for $k \in I$, and $\psi(\{\phi_i(t)\}_i)
= t$. The former follows from
\[\phi_k(\psi(\{s_i\}_i)) \ =\ \phi_k(\alpha\{i \mapsto s_i\}_i)
\ =\ \alpha[k]1_{X_k} ;_\alpha \alpha\{i \mapsto s_i\}_i
\ =\ 1_{X_k} ;_\alpha s_k \ =\ s_k
\]
and the latter from
\[\psi(\{\phi_i(t)\}_i) \ =\ \psi(\alpha[i]1_{X_i} ;_\alpha t)
\ =\ \beta\{i \mapsto \alpha[i]1_{X_i} ;_\alpha t\}_i
\ =\ \beta\{i \mapsto \alpha[i]1_{X_i} \}_i ;_\alpha t
\ =\ 1_{\Sigma_i X_i} ;_\alpha t \ =\ t
\]

It remains to show that this correspondence is natural, however, this
follows immediately from the rewrites (3) and (4).
\end{proof}

\subsection{Representability} \label{sec-repre}

A multi-map $\Gamma \to^u X$ is said to \textbf{represent $\Gamma$ as input}
(cf.~\cite{hermida00:representable,cockett03:polybicat}), if cutting with $u$
at $X$ induces a natural bijection as follows:

\[\infer={\Gamma_1,X,\Gamma_2 \ra \Delta}{\Gamma_1,\Gamma,\Gamma_2 \ra \Delta}
\]
Dually, a comulti-map $Y \to^v \Gamma$ \textbf{represents $\Gamma$ as
output}, if cutting with $v$ at $Y$ induces an analogous bijection. A
polycategory is called \textbf{representable} if each sequence of formulas
is representable as input and as output.

\begin{theorem}
$\channel(\cat{A})$ is representable.
\end{theorem}

\begin{proof}
We will prove that a bundle of channels $\{\beta_i \cc X_i\}_i$ is
represented as input by the following multi-map

\senta{$\alpha\left\<\alpha_i \mid \beta_i \mapsto \beta_i \equiv_{X_i}
\alpha_i \right\>::\{\beta_i \cc X_i\}_i \ra \alpha \cc \bigox_i \alpha_i
\cc X_i$}
The symmetry of the term calculus will then deliver representability.

This means that we must show the correspondences
\[\vcenter{\infer={\Gamma,\alpha \cc \bigox_i \alpha_i \cc X_i \ra \Delta}
{\Gamma,\{\alpha_i \cc X_i\}_i \ra \Delta}} \qquad \text{and} \qquad
\vcenter{\infer={\Gamma,\alpha \cc \top \ra \Delta}{\Gamma \ra \Delta}}
\]
are bijective and natural. The proof of the non-nullary case follows in a
manner very similar to that of Proposition~\ref{prop-hassums}.

Suppose we have maps $t::\Gamma,\alpha \cc \bigox_i \alpha_i \cc X_i \ra
\Delta$ and $s::\Gamma,\{\alpha_i \cc X_i\}_i \ra \Delta$. The map $s$ can,
via the ltensor derivation, be used to construct a map

\senta{$\psi(s) = \alpha\<(\alpha_i)_i \mapsto s\>::
\Gamma,\alpha \cc \bigox_i \alpha_i \cc X_i \ra \Delta$}
and $t$ may be cut with the representing multi-map to get a map

\senta{$\phi(t) = \alpha\left\<\alpha_i \mid \beta_i \mapsto \beta_i
\equiv_{X_i} \alpha_i \right\>_i ;_\alpha t ::
\Gamma,\{\alpha_i \cc X_i\}_i \ra \Delta$}

To prove that this correspondence is bijective we must establish both
$\phi(\psi(s)) = s$ and $\psi(\phi(t)) = t$. Suppose $I = \{1,\ldots,n\}$.
The former follows from
\begin{align*}
\phi(\psi(s)) &= \phi(\alpha\<(\alpha_i)_i \mapsto s\>) \\
&= \alpha\<\alpha_i \mid \beta_i \mapsto \beta_i \equiv_{X_i} \alpha_i\>_i
;_\alpha \alpha\<(\alpha_i)_i \mapsto s\>) \\
&= \beta_n \equiv_{X_n} \alpha_n ;_{\alpha_n}( \cdots (\beta_1 \equiv_{X_1}
\alpha_1 ;_{\alpha_1} s) \cdots ) \\
&= s
\end{align*}
and the latter from
\begin{align*}
\psi(\phi(t)) &= \psi(\beta\<\beta_i \mid \alpha_i \mapsto \alpha_i
\equiv_{X_i} \beta_i\>_i ;_\beta t) \\
&= \alpha\<(\alpha_i)_i \mapsto \beta\<\beta_i \mid \alpha_i \mapsto
\alpha_i \equiv_{X_i} \beta_i\>_i ;_\beta t\> \\
&= \alpha\<(\alpha_i)_i \mapsto \beta\<\beta_i \mid \alpha_i \mapsto
\alpha_i \equiv_{X_i} \beta_i\>_i\> ;_\beta t \\
&= 1_{\bigox_i X_i} ;_\beta t \\
&= t
\end{align*}

To establish the nullary case suppose we have maps $t::\Gamma,\alpha \cc
\top \ra \Delta$ and $s::\Gamma \ra \Delta$, and define
\[\psi(s) = \alpha\<(\,) \mapsto s\> \quad \text{and} \quad
\phi(t) = \alpha\<\,\> ;_\alpha t
\]
To see that it is again a bijection observe
\[\phi(\psi(s)) \ =\ \phi(\alpha\<(\,) \mapsto s\>)
\ =\ \alpha\<\,\> ;_\alpha \alpha\<(\,) \mapsto s\> \ =\ s
\]
and 
\[\psi(\phi(t)) \ =\ \psi(\beta\<\,\> ;_\beta t)
\ =\ \alpha\<(\,) \mapsto \beta\<\,\> ;_\beta t\>
\ =\ \alpha\<(\,) \mapsto \beta\<\,\>\> ;_\beta t
\ =\ 1_{\top} ;_\beta t \ =\ t
\]

It is left to show that the correspondence is natural, however, this follows
immediately from the rewrites (7) and (8).
\end{proof}

\subsection{The free additive linearly distributive category}
\label{sec-free}

The goal of this section is to prove:

\begin{theorem}
$\channel(\cat{A})$ is the free representable polycategory with additives generated from \cat{A}.
\end{theorem}

As any linearly distributive category with additives generates a
polycategory with additives, an immediate corollary of this is:

\begin{corollary}
The maps of $\channel(\cat{A})$ form the free linearly distributive category
with additives.
\end{corollary} 

In order to prove this theorem we must show that in any representable
polycategory with additives the identities (1)-(24) must hold.

It is routine to show (see \cite{pastro:msc}) that a polycategory has
poly-sums if and only if it has a cotupling operation which is distributive
on the non-sum channel, and injections which when cut against a cotuple on
the sum channel delivers the appropriate component of the cotuple, and
finally it satisfies ``surjective pairing'' which is the requirement that
the cotuple of the injections is the identity.  

First observe that the injection term $\alpha[a]f$ can be translated as
$f;b_k$ where $b_k::X_k \ra \sum_i X_i$, for $k \in I$. This is a valid
identification as there is a reduction of derivations:

\[\vcenter{\infer{\Gamma \vd \sum_i X_i, \Delta}{\Gamma \vd X_k,\Delta &
\infer{X_k \vd \sum_i X_k}{X_k \vd X_k}}} \quad \Lra \quad
\vcenter{\infer{\Gamma \vd \sum_i X_i, \Delta}{\Gamma \vd X_k,\Delta}}
\]
The term calculus identities for the additives express precisely the
requirements described above with the exception of surjective pairing.
However, recall that this is implicit in the term calculus from the manner
in which the identity map for the coproduct is defined.

The only remaining difficulty is to translate the ``forking'' construct of
the term logic. For this suppose $I = \{1,\ldots,n\}$ for $n \geq 1$, then
in any representable polycategory we may translate the map $\alpha\<\alpha_i
\mid \Lambda_i \mapsto f_i\>_i$ as $f_n ;_{\alpha_n} (\ldots
(f_1 ;_{\alpha_2} r) \ldots)$ where $r$ is the representing multi-map. 

\[\infer{\Gamma \vd \bigox_i X_i,\Delta}{
    \infer{\Gamma_1 \vd X_1,\Delta_1}{f_1} &
\infer{X_1,\Gamma_2,\cdots,\Gamma_n \vd \bigox_i X_i,\Delta_2,\cdots,\Delta_n}{
\infer{\ \cdots \ }{
    \infer{\Gamma_n \vd X_n,\Delta_n}{f_n} &
    \infer={\{X_i\}_i \vd_r \bigox_i X_i}{
       \infer{\bigox_i X_i \vd \bigox_i X_i}{1_{\bigox_i X_i}}}}}}
\]

The identities are now an easy consequence. In particular, (9), (10), (13),
and (14) are a direct consequence of the translation of ``forking'' described 
above.

\section{The decision procedure} \label{sec-dp}

Cut elimination seen as a term rewriting supplies a way of rewriting the
terms (in the initial calculus) to remove the cut completely. In
$\channel{(\cat{A})}$ it allows the cuts to be moved into the underlying
polycategory $\cat{A}$. This relies on the following routine but nonetheless
technical observation which also delivers a ``categorical cut elimination''
theorem as the equality of proofs (as determined by the equalities above) is
now guaranteed to be maintained by the elimination procedure.

\begin{theorem} In $\channel{(\cat{A})}$:
\begin{enumerate}[{\upshape (i)}]
\item The rewriting on maps given by (1)-(14) terminates.
\item The rewriting on maps given by (1)-(14) is confluent modulo the
permuting conversions (15)-(24).
\end{enumerate}
\end{theorem}

The proof is presented in some detail in Appendix~\ref{sec-proof-ce}. The
argument involves resolving all the critical pairs (rewriting against
rewriting and rewriting against permuting conversion) modulo the permuting
conversions.  To show that the rewriting terminates in forms which are
equivalent with respect to the permuting conversions involves showing that
the resolution of each critical pair always reduces the cost of the
frontier. This does involve some subtlety as the nullary additive permuting
conversions can produce some apparently non-reducing steps. However,
somewhat surprisingly, a very basic multiset measure of cut heights applied
carefully does suffice for the whole argument to go through quite smoothly. 

The effect of this theorem is to provide a decision procedure modulo the
decidability of the underlying polycategory $\cat{A}$. The rewriting
normalizes terms by moving the cut into the non-logical axioms. The
equivalence of terms is then determined by the decision procedure in
$\cat{A}$ and the permuting conversions. The subformula property delivered
by cut elimination ensures that there are only finitely many proofs which 
do not involve cut (excluding how the non-logical steps are filled). This
means that, in principal, in order to decide the equality of two terms one
can simply search all the terms equivalent by permuting conversions to one of
the proofs for a term equivalent to the other. Two terms are equivalent
when the proof structure of the two terms match exactly and the non-logical
steps which correspond are equivalent. This observation is sufficient to
establish the following formal result:

\begin{corollary}
The equivalence of proofs in $\channel{(\cat{A})}$ is decidable whenever the equivalence of 
proofs in $\cat{A}$ is decidable.
\end{corollary}

Recall that all discrete polycategories and free polycategories (including
those with negation) are decidable so that there is a ready source of
decidable polycategories.

Of course, the procedure proposed is highly inefficient and we wish to finish
the paper with some remarks on how one can make it more efficient. In order
to decide the equality of two terms it is worth putting one term aside as a
template. However, it is also sensible to ensure that this term has all
tensor/par splitting done first and then all the tupling/cotupling. This
structure can always be pulled up to the root of the term whenever an active
multiplicative or additive type is present.  This leaves the tensor/par
forking and the injections/projections. The latter can be pulled up in a
similar manner to \cite{cockett01:finite}. Pulling up a forking to match a
template term which starts with such is more complex: this is exactly where
the complexity of rewirings (due to the multiplicative units) in the purely
multiplicative fragment bites. The complexity of solving this step remains
an open question.
 


\clearpage
\appendix
\section{The units}

The way in which the cut elimination procedure handles the reductions and
permuting conversions when the index sets are empty can be quite subtle.
In the case of ltensor or rpar, however, it is very straightforward,
however the other cases are quite tricky. To clarify this, in this section
we make these special cases explicit.

The annotated nullary versions of the inference rules are:

\begin{center}
\ovalbox{\parbox{80ex}{
\[\begin{array}{ccc}
\infer[\text{(cotuple)}]{\alpha\{\,\} ::\Gamma,\alpha \cc \mathbf{0} \vd
\Delta}{}
&& \infer[\text{(tuple)}]{\alpha\{\,\}::\Gamma \vd \alpha \cc \mathbf{1},
\Delta}{}
\medskip\\
\infer[\text{(ltensor)}]
{\alpha\<(\,) \mapsto f\>::\Gamma,\alpha \cc \top \vd \Delta}
{f::\Gamma \vd \Delta} &&
\infer[\text{(rpar)}]
{\alpha\<(\,) \mapsto f\>::\Gamma \vd \alpha \cc \bot,\Delta}
{f::\Gamma \vd \Delta}
\medskip\\
\infer[\text{(lpar)}]{\alpha\<\,\>:: \alpha \cc \bot \vd\ }{} &&
\infer[\text{(rtensor)}]{\alpha\<\,\>::\ \vd \alpha \cc \top}{}
\end{array}
\]
}}\end{center}

In the cotuple and typle rules the notation is ambiguous as one cannot
derive the context from the terms. To correct this we shall write the terms
above as $\alpha\{\,\}_{\Gamma \vd \Delta}$. An imporant observation is that
the nullary rtensor and nullary lpar rules may only be applied if there are
no other ``active'' channels present. Considering this, there are three
(non-dual) reductions that are relevant to this setting corresponding to the
rewrites (3), (7), and (13). The reduction (9) is not relavant to this setting
as it will reduced to (13). Given a term $f:: \gamma \cc X,\Delta_1 \vd
\Delta_2$ the reductions are as follows:
\[\begin{array}{rrcl}
(3) & \alpha\{\,\}_{\Gamma_1 \vd \Gamma_2,\gamma \cc X}\ ;_\gamma f
&\Lra& \alpha\{\,\}_{\Gamma_1,\Delta_1 \vd \Gamma_2,\Delta_2} \medskip\\
(7) & \alpha\<(\,) \mapsto g\>\ ;_\gamma f &\Lra&
\alpha\<(\,) \mapsto g ;_\gamma f\> \medskip\\
(13) & \gamma\<\,\>\ ;_\gamma \gamma\<(\,) \mapsto f\> &\Lra& f
\end{array}
\]

With the interchange rules the typing of the term will determine which
interchanges are allowed to take place. For example, consider (23), the
ltensor-rtensor identity:
\[\alpha\left\<(\alpha_i)_i \mapsto \beta\left\<\begin{array}{l}
\beta_j | \Lambda_j \mapsto g_j \\
\beta_k | \Lambda_k \cup (\bigcup_i \{\alpha_i\}) \mapsto f
\end{array}\right\>_{j \neq k}\right\>
\pc
\beta\left\<\begin{array}{l}
\beta_j | \Lambda_j \mapsto g_j \\
\beta_k | \Lambda_k \cup \{\alpha\} \mapsto \alpha\<(\alpha_i)_i \mapsto f\>
\end{array}\right\>_{j \neq k}
\]
On the left-hand side the only time the nullary rtensor rule may be applied
is if $\bigcup_i \{\alpha_i\} = \emptyset$ and no other channels are ``active''.
This corresponds to the term
\[\alpha\<(\,) \mapsto \beta\<\,\>\>:: \alpha \cc \top \vd \beta \cc \top
\]
However, on the right-hand side the nullary rtensor rule will never apply
as the channel $\alpha$ is still ``active''.  Thus, in this case, these
rules may not be interchanged. If only the ltensor is empty (on channel
$\alpha$) then the rules may freely be interchanged.

Considering the typing of the terms this leaves 13 permuting conversion:
three versions of (15) (corresponding to the cases when only $I = \emptyset$,
only $J = \emptyset$, and both $I = J = \emptyset$), (16), three variants
of (17), (18) (when $I = \emptyset$), (20), three variants of (22), and
finally (23) (when $I = \emptyset$). Fortunately, the form of the rewrites
are quite similar, and so we present only the three variants of (15), (20),
and (23). For simplicity we assume that $\alpha$ is a domain channel and
$\beta$ is a codomain channel. Here we drop the typing on the term and
indicate to the right in brackets.
\[\begin{array}{rrcll}
(15) & \alpha\{a_i \mapsto \beta\{\,\}\}_i & \pc & \beta\{\,\} &
(\Gamma,\alpha \cc \sum X_i \vd \beta \cc \mathbf{1},\Delta) \medskip\\
(15) & \alpha\{\,\} & \pc & \beta\{b_j \mapsto \alpha\{\,\}\}_j &
(\Gamma,\alpha \cc \mathbf{0} \vd \beta \cc \prod Y_j,\Delta) \medskip\\
(15) & \alpha\{\,\} & \pc & \beta\{\,\} &
(\Gamma,\alpha \cc \mathbf{0} \vd \beta \cc \mathbf{1},\Delta) \medskip\\
(20) & \alpha[a]\beta\<(\,) \mapsto f\> &\pc& \beta\<(\,) \mapsto \alpha[a]f\>&
(\Gamma,\alpha \cc \prod X_i \vd \beta \cc \bot,\Delta) \smallskip\\
(23) & \scriptstyle{\alpha\left\<() \mapsto \beta\left\<\begin{array}{l}
\scriptstyle{\beta_j | \Lambda_j \mapsto g_j} \\
\scriptstyle{\beta_k | \Lambda_k \mapsto f}
\end{array}\right\>_{j \neq k}\right\>}
&\pc&
\scriptstyle{\beta\left\<\begin{array}{l}
\scriptstyle{\beta_j | \Lambda_j \mapsto g_j} \\
\scriptstyle{\beta_k | \Lambda_k \cup \{\alpha\} \mapsto \alpha\<() \mapsto f\>}
\end{array}\right\>_{j \neq k}} &
(\Gamma,\alpha \cc \top \vd \beta \cc \bigox Y_j,\Delta)
\end{array}
\]

\section{Proof of cut elmination and the Church-Rosser property}
\label{sec-proof-ce}

In this appendix we show that the cut elimination procedure terminates and
that the rewrite system induced by the cut elimination rewrites and the
communting conversion has the Church-Rosser property. We begin by proving
that cut elimination is a terminating procedure.

\begin{remark}[Generalized rewrites] \label{remark-genre}
The number of rewritings in our system motivates the use of a ``generalized''
system of rewritings which help to reduce this number. Let $\alpha(f)$
denote any of the morphims
\[\alpha\{a_i \mapsto f_i\}_i, \qquad \alpha[a]f, \qquad
\alpha\<(\alpha_i) \mapsto f\>, \qquad
\alpha\<\alpha_i \mid \Lambda_i \mapsto f_i\>_i
\]
Then $\alpha(f);_\gamma g \Lra \alpha(f ;_\gamma g)$ will be used
respectively to denote:
\[\begin{array}{rcl}
\alpha\{a_i \mapsto f_i\}_i ;_\gamma g &\xymatrix{\ar@{=>}[r]^{(3)}&}&
  \alpha\{a_i \mapsto f_i ;_\gamma g\}_i \medskip\\
\alpha[a]f ;_\gamma g &\xymatrix{\ar@{=>}[r]^{(5)}&}&
  \alpha[a](f ;_\gamma g) \medskip\\
\alpha\<(\alpha_i) \mapsto f\> ;_\gamma g &\xymatrix{\ar@{=>}[r]^{(7)}&}&
  \alpha\<(\alpha_i) \mapsto f;_\gamma g\> \smallskip\\
\alpha\<\alpha_i \mid \Lambda_i \mapsto f_i\>_i ;_\gamma g
&\xymatrix{\ar@{=>}[r]^{(9)}&}&
\alpha\left\<\begin{array}{l}
\alpha_i \mid \Lambda_i \mapsto f_i \\
\alpha_k \mid \Lambda_k \cup \Lambda(g) \setminus \{\gamma\}
\mapsto f_k ;_\gamma g
\end{array}\right\>_{i \neq k}
\end{array}
\]

Dually $\xymatrix@1{g ;_\gamma \alpha(f) \ar@{=>}[r] & \alpha(g;_\gamma f)}$
will denote any of the rewrites (4), (6), (7), or (10). The rewrites (11)
and (13) (and their duals) may be represented as $\xymatrix@1{\gamma(f)
;_\gamma \gamma(g) \ar@{=>}[r] & f ;_\gamma g}$ and the permuting
conversions (15) through (24) may be represented as
$\xymatrix@1{\alpha(\beta(f)) \ar@{|=|}[r] & \beta(\alpha(f))}$.
\end{remark}

\subsection{The cut measure on terms}

The purpose of this section is to show that the cut elimination procedure
terminates. To this end we define a bag of cut heights and show that the
bag is strictly reduced on each of the cut elmination rewrites.

We begin by defining the multiset ordering of Dershowitz and
Manna~\cite{dershowitz79:proving}. Let $(S,\succ)$ be a partially-ordered
set, and let $\script{M}(S)$ denote the multisets (or bags) over $S$. For
$M,N \in \script{M}(S)$, $M > N$ (``$>$'' is called the \textbf{multiset}
(or \textbf{bag}) \textbf{ordering}), if there are multisets $X,Y \in
\script{M}(S)$, where $\emptyset \neq X \subseteq M$, such that
\[N = (M \bs X) \cup Y \quad \text{and} \quad (\forall y \in Y)
(\exists x \in X)\ x \succ y
\]
where $\cup$ here is the multiset union.

For example,
\[[3] > [2,2,2,1], \quad [7,3] > [7], \quad [5,2] > [5,1]
\]

Recall from~\cite{dershowitz79:proving} that if $(S,\succ)$ is a total order
(linear order) then $\script{M}(S)$ is a total order. To see this consider
$M,N \in \script{M}(S)$. To determine whether $M > N$ sort the elements of
both $M$ and $N$ and then compare the two sorted sequences lexicograpically.

We now define the \textbf{height} of a term as:

\begin{itemize}
\item $\hgt[a] = 1$ when $a$ is an atomic map (or an identity)
\item $\hgt[\alpha\{a_i \mapsto f_i\}_{i \in I}] = 1+ \max\{\hgt[f_i] \mid
i \in I\}$
\item $\hgt[\alpha[a_k] \cdot f] = 1+\hgt[f]$
\item $\hgt[\alpha\<(\alpha_i)_{i \in I} \mapsto f\>] = 1 + \hgt[f]$
\item $\hgt[\alpha\<\alpha_i \mid \Omega_i \mapsto f_i\>_{i \in I}] =
1 + \sum_{i \in I} \hgt[f_i]$
\item $\hgt[f ; g] = \hgt[f] + \hgt[g]$
\end{itemize}

The \textbf{height of a cut} is defined simply as its height, e.g.,
$\cuthgt[f;g] = \hgt[f;g]$. Define a function $\Lambda:T \ra \bag(\mathbb{N})$
which takes a term to its bag of cut heights.

\begin{proposition} \quad
\begin{enumerate}[{\upshape (i)}]
\item If $\xymatrix{t_1 \ar@{=>}[r]& t_2}$ then $\Lambda(t_1) > \Lambda(t_2)$.
\item If $\xymatrix{t_1 \ar@{|=|}[r]^{(a)}& t_2}$ and $(a)$ is an
interchange which does not involve the nullary cotuple or tuple then
$\Lambda(t_1) = \Lambda(t_2)$.
\end{enumerate}
\end{proposition}

\begin{proof}
We begin with the proof of part (i). There are three properties that must be
shown: $\hgt[t_1] \geq \hgt[t_2]$, the height of each non-principal cut does
not increase, and the height of any cut produced from the principal cut is
strictly less than the height of the principal cut.

A simple examination of the rewrites will confirm
that if $t_1 \Lra t_2$ then $\hgt[t_1] \geq \hgt[t_2]$:

\begin{equation*} \tag*{(1), dually (2)}
\hgt[f;1] = \hgt[f] + \hgt[1] > \hgt[f]
\end{equation*}

\begin{align*} \tag*{(3), dually (4)}
\hgt[\alpha\{a_i \mapsto f_i\}_{i \in I} ;_\gamma g]
&= \hgt[\alpha\{a_i \mapsto f_i\}_{i \in I}] + \hgt[g] \\
&= 1+ \max\{\hgt[f_i] \mid i \in I\} + \hgt[g] \\
&= 1+ \max\{\hgt[f_i] + \hgt[g] \mid i \in I\}  \\
&= \hgt[\alpha\{a_i \mapsto f_i ;_\gamma g\}_{i \in I}]
\end{align*}

If $I = \emptyset$ then 
\[\hgt[\alpha\{\,\} ;_\gamma g] =
\hgt[\alpha\{\,\}] + \hgt[g] > \hgt[\alpha\{\,\}]
\]

\begin{align*} \tag*{(5), dually (6)}
\hgt[\alpha[a_k]f ;_\gamma g] &= \hgt[\alpha[a_k]f]+\hgt[g] \\
&= 1+ \hgt[f] + \hgt[g] \\
&= \hgt[\alpha[a_k](f ;_\gamma g)]
\end{align*}

\begin{align*} \tag*{(7), dually (8)}
\hgt[\alpha\<(\alpha_i)_{i \in I} \mapsto f\> ;_\gamma g]
&= \hgt[\alpha\<(\alpha_i)_{i \in I} \mapsto f\> +\hgt[g]  \\
&= 1 + \hgt[f] + \hgt[g] \\
&= \hgt[\alpha\<(\alpha_i)_{i \in I} \mapsto f ;_\gamma g\>
\end{align*}

\begin{align*} \tag*{(9), dually (10)}
\hgt[\alpha\<\alpha_i \mid \Omega_i \mapsto f_i\>_{i \in I} ;_\gamma g]
&= \hgt[\alpha\<\alpha_i \mid \Omega_i \mapsto f_i\>_{i \in I}]+\hgt[g]\\
&= 1 + \sum_{k \neq i \in I} \hgt[f_i] + \hgt[f_k] + \hgt[g] \\
&= \hgt\left[\alpha\left\<\begin{array}{l}
\alpha_i \mid \Omega_i \mapsto f_i \\
\alpha_k \mid \Omega_k \mapsto f_k ;_\gamma g
\end{array} \right\>_{k \neq i \in I}\right]
\end{align*}

\begin{align*} \tag*{(11), dually (12)}
\hgt[\gamma[a_k]f ;_\gamma \gamma\{a_i \mapsto g_i\}_{i \in I}]
&= \hgt[\gamma[a_k]f] + \hgt[\gamma\{a_i \mapsto g_i\}_{i \in I}] \\
&= 1+ \hgt[f] + 1 + \max\{\hgt[g_i] \mid i \in I\} \\
&> \hgt[f] + \max\{\hgt[g_i] \mid i \in I\} \\
&\geq \hgt[f] + \hgt[g_k] \\
&= \hgt[f ;_\gamma g_k]
\end{align*}

\begin{align*} \tag*{(13), dually (14)}
\hgt[\gamma\<\alpha_i \mid \Omega_i \mapsto f_i\>_i ;_\gamma
\gamma\<(\alpha_i)_i \mapsto g\>]
&= \hgt[\gamma\<\alpha_i \mid \Omega_i \mapsto f_i\>_i] +
\hgt[\alpha\<(\alpha_i)_i \mapsto g\>] \\
&= 1 + \sum_i \hgt[f_i] + 1 + \hgt[g] \\
&> \sum_i \hgt[f_i] + \hgt[g] \\
&= \hgt\left[f_n ;_{\alpha_n} ( \cdots (f_2 ;_{\alpha_2} (f_1 ;_{\alpha_1} g))
\cdots)\right]
\end{align*}

If $I = \emptyset$ then
\begin{align*}
\hgt[\gamma\<\,\>_i ;_\gamma \gamma\<(\,) \mapsto g\>]
&= \hgt[\gamma\<\,\>_i] + \hgt[\gamma\<(\,) \mapsto g\>] \\
&= 1+1+\hgt[g] \\
&> \hgt[g]
\end{align*}

Moreover, this implies that cuts below and cuts above the redex will not
increase their cut height on a rewriting. 

Finally, consider the principal cut of the reduction. Rewrite (1) (dually (2))
removes a cut and so strictly reduces the bag of cut heights. It is an easy
observation that (5), (7), (9), and (11) (and their duals) each replace a
cut with one of lesser height, and that (3) and (13) (and their duals)
replace a cut with zero or more cuts of lesser height. Thus applying any of
the rewrites strictly reduces the bag.

We now prove part (ii). For the equations (15), (16), (17), and (18) we
assume that the index sets are non-empty. This then implies that the
commuting conversions are all of the form $\alpha(\beta(f))$ and thus
\[\hgt[\alpha(\beta(f))] = 1 + \hgt[\beta(f)] = 1 + 1+ \hgt[f] = 
\hgt[\beta(\alpha(f))]
\]
which proves that the height does not change across these (non-empty
(co)tuple) interchanges.
\end{proof}

To see that the height is not invariant across the emtpy cotuple (dually the
tuple) rule recall one of the nullary versions of the rewrite (15):
\[\alpha\{\,\} \pc \beta\{b_j \mapsto \alpha\{\,\}\}_j
\]
The height on the left-hand side is one, while on the right-hand side the
height is two.

\subsection{Proof of the Church-Rosser property} \label{sec-church-rosser}

In this section we present a proof of the Church-Rosser property for
morphisms. We wish to show that given any two morphisms related by a series
of reductions and permuting conversions
\[\xymatrix{t_1 \ar@{<=}[r] & t_2 \ar@{|=|}[r] & t_3 \ar@{=>}[r] &
\quad \cdots \quad &
\ar@{=>}[l] t_{n-2} \ar@{|=|}[r] & t_{n-1} \ar@{=>}[r] & t_n}
\]
there is an alternative way of arranging the reductions and permuting
conversions so that $t_1$ and $t_n$ can be reduced to terms which are
related by the permuting conversions alone. That is, we wish to show that
there is a convergence of the following form:
\[\xymatrix{t_1 \ar@{=>}[dr]_{*} & && & t_n \ar@{=>}[dl]^{*} \\
& t_1' \ar@{|=|}[rr]_{*} && t_n'}
\]

When the rewriting system terminates (in the appropriate sense) this allows
the decision procedure for the equality of $\Sigma\Pi$-terms to be reduced
to the decision procedure for the permuting conversions (see
Section~\ref{sec-dp}). In order to test the equality of two terms, one can
rewrite both terms into a reduced form (one from which there are no further
reductions), and these will be equal if and only if the two reduced forms
are equivalent through the permuting conversions alone. In the current
situation the reduction process is, of course, the cut-elimination procedure.

Following~\cite{cockett01:finite} we say a rewrite system is \textbf{locally
confluent modulo equations} if any (one step) divergence of the following
form
\[\vcenter{\xymatrix{& t_0 \ar@{=>}[dl] \ar@{=>}[dr] \\ t_1 && t_2}}
\qquad \text{or} \qquad
\vcenter{\xymatrix{& t_0 \ar@{=>}[dl] \ar@{|=|}[dr] \\ t_1 && t_2}}
\]
(where ``$\xymatrix@1{\ar@{=>}[r]&}$'' denotes a reduction and
``$\xymatrix@1{t_1 \ar@{|=|}[r]& t_2}$'' an equation) has a convergence,
respectively, of the form
\[
\vcenter{\xymatrix{t_1 \ar@{|=>}[dr]_{*} && t_2 \ar@{|=>}[dl]^{*} \\ & t'}}
\qquad \text{and} \qquad
\vcenter{\xymatrix{t_1 \ar@{|=>}[ddr]_{*} && t_2 \ar@{=>}[d]  \\
&& t_2' \ar@{|=>}[dl]^{*} \\ & t'}}
\]
where the new arrow ``$\xymatrix@1{\ar@{|=>}[r]&}$'' indicates either
an equality or a reduction in the indicated direction.

This gives:
\begin{proposition} \label{prop-conf}
Suppose $(N,\mathcal{R},\mathcal{E})$ is a rewriting system with the
equations equipped with a well-ordered measure on the rewrite arrows such
that the measure of the divergences is strictly greater than the measure of
the convergences then the system is confluent modulo equations if and only if
it is locally confluent modulo equations.
\end{proposition}

\begin{proof}
If the system is confluent modulo equations it is certainly locally
confluent modulo equations. Conversely suppose we have a chain of reductions,
equations, and expansions. We may associate with it the bag of measures of
the arrows of the sequence.

The idea will be to show that replacing any local divergence in this chain
by a local confluence will result in a new chain whose bag measure is
strictly smaller. However, this can be seen by inspection as we are removing
the arrows associated with the divergence and replacing them with the arrows
associated with the convergence. The measure on the arrows associated with
the divergence is strictly greater then that of the measure on the arrows
associated with the convergence.

Thus, each rewriting reduces the measure and, therefore, any sequence
of rewriting on such a chain must terminate. However, it can only terminate
when there are no local divergences to resolve. This then implies that the
end result must be a confluence modulo equations.
\end{proof}

\subsubsection{Resolving critical pairs locally}

The proof of Church-Rosser involves examining all the possible critical 
pairs involving reductions or reductions and conversions, and showing
that they are all of the form shown above and that they may be resolved in
the way shown above. It then must be shown that there is some measure on the
arrows which decreases when replacing a divergence with a convergences.
This will then suffice to show that our system is locally confluent modulo
equations, so that by Proposition~\ref{prop-conf}, it is confluent modulo
equations. The rewrites (1)-(12) are the ``reductions'' and the commuting
conversions (13)-(24) are the ``equations''. The resolutions of the critical
pairs will be presented using the ``generalized'' rewrites (see
Remark~\ref{remark-genre}). For the additive rewrites see~\cite{pastro:msc}
where they have been written out in detail.

The reductions are as follows. Note that these reductions assume that all
index sets are non-empty. The reductions when the index set is empty are
handled separately below.

\begin{description}
\item[Reduction diagram 0:] $\xymatrix{1;1 \ar@{=>}[r]_{(2)}^{(1)} & 1}$

\item[Reduction diagram 1:] Substitute (3), (5), (7), or (9) for $(a)$ to
get the reduction diagrams for (1)-(3), (1)-(5), (1)-(7), and (1)-(9).
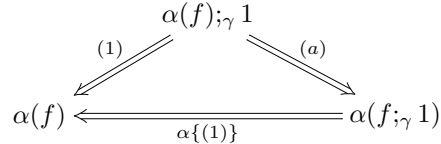
\[\xymatrix@R=5ex@C=7ex{
& \alpha(f) ;_\gamma 1 \ar@{=>}[dl]_-{(1)} \ar@{=>}[dr]^-{(a)} \\
\alpha(f) && \alpha(f ;_\gamma 1) \ar@{=>}[ll]^-{\alpha\{(1)\}} }
\]
Substituting the dual rewrites in the mirror image of the diagram above give
the dual reduction diagrams.

\item[Reduction diagram 2:] each row in the table corresponds to the
resolution of the critical pair ($a$)-($b$).
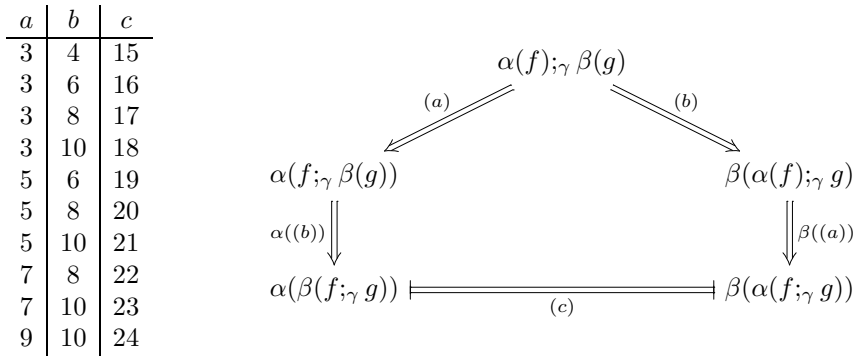
\[\begin{array}{c|c|c}
a & b & c \\ \hline
3 & 4 & 15 \\
3 & 6 & 16 \\
3 & 8 & 17 \\
3 & 10 & 18 \\
5 & 6 & 19 \\
5 & 8 & 20 \\
5 & 10 & 21 \\
7 & 8 & 22 \\
7 & 10 & 23 \\
9 & 10 & 24
\end{array}
\qqqquad
\vcenter{\xymatrix@M=1ex@C=20ex@R=10ex@!0{
& \alpha(f) ;_\gamma \beta(g) \ar@{=>}[dl]_{(a)} \ar@{=>}[dr]^{(b)} \\
\alpha(f ;_\gamma \beta(g)) \ar@{=>}[d]_{\alpha((b))}
&& \beta(\alpha(f) ;_\gamma g) \ar@{=>}[d]^{\beta((a))} \\
\alpha(\beta(f ;_\gamma g)) \ar@{|=|}[rr]_{(c)}
&& \beta(\alpha(f ;_\gamma g))}}
\]

\item[Reduction diagram 3:] each row in the table corresponds to the
resolution of the critical pair ($a$)-($b$).
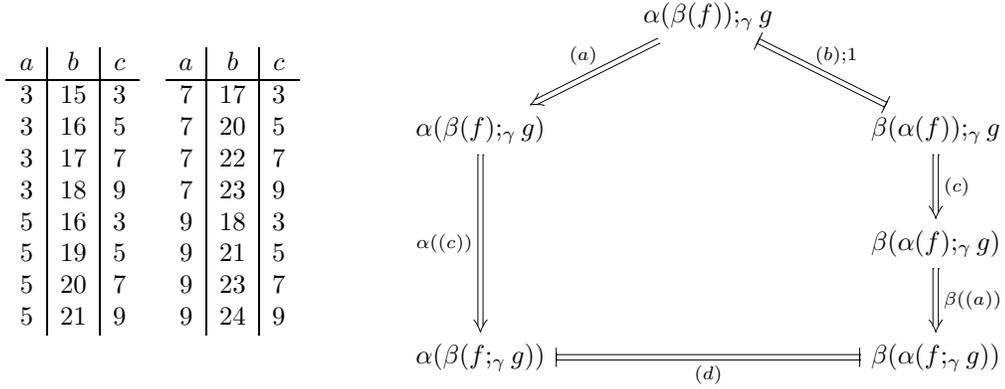
\[\begin{array}{c|c|c}
a & b & c \\ \hline
3 & 15 & 3 \\
3 & 16 & 5 \\
3 & 17 & 7 \\
3 & 18 & 9 \\
5 & 16 & 3 \\
5 & 19 & 5 \\
5 & 20 & 7 \\
5 & 21 & 9 \\
\end{array}
\quad
\begin{array}{c|c|c}
a & b & c \\ \hline
7 & 17 & 3 \\
7 & 20 & 5 \\
7 & 22 & 7 \\
7 & 23 & 9 \\
9 & 18 & 3 \\
9 & 21 & 5 \\
9 & 23 & 7 \\
9 & 24 & 9
\end{array}
\qqqquad
\vcenter{\xymatrix@M=1ex@C=20ex@R=10ex@!0{
& \alpha(\beta(f)) ;_\gamma g \ar@{=>}[dl]_{(a)} \ar@{|=|}[dr]^{(b);1} \\
\alpha(\beta(f) ;_\gamma g) \ar@{=>}[dd]_{\alpha((c))}
&& \beta(\alpha(f)) ;_\gamma g \ar@{=>}[d]^{(c)} \\
&& \beta(\alpha(f) ;_\gamma g) \ar@{=>}[d]^{\beta((a))} \\
\alpha(\beta(f ;_\gamma g)) \ar@{|=|}[rr]_{(d)}
&& \beta(\alpha(f ;_\gamma g))}}
\]

\item[Reduction diagram 4:] each row in the table corresponds to the
resolution of the critical pair ($a$)-($b$).
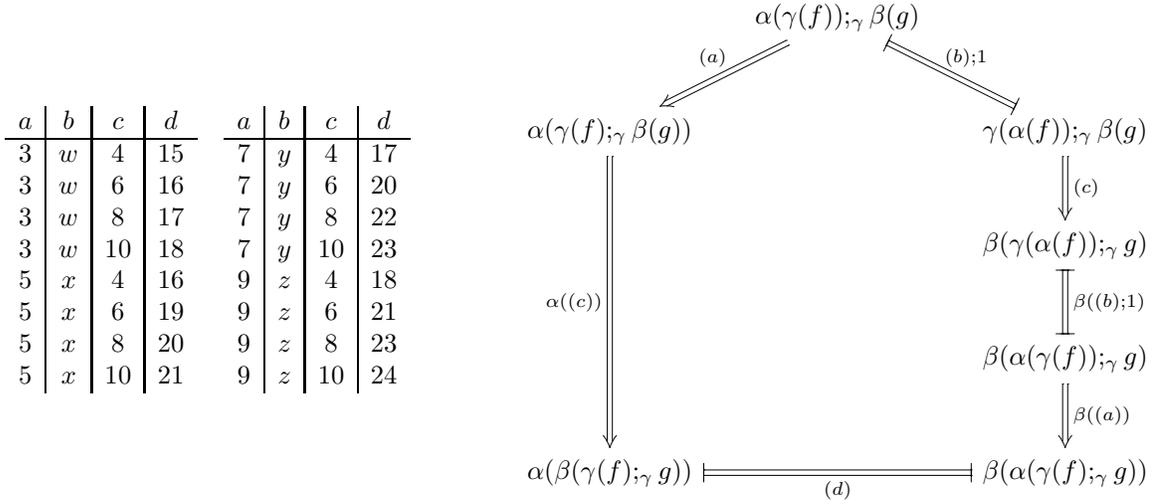
\[\begin{array}{c|c|c|c}
a & b & c & d\\ \hline
3 & w & 4 & 15 \\
3 & w & 6 & 16 \\
3 & w & 8 & 17 \\
3 & w & 10 & 18 \\
5 & x & 4 & 16 \\
5 & x & 6 & 19 \\
5 & x & 8 & 20\\
5 & x & 10 & 21
\end{array}
\quad
\begin{array}{c|c|c|c}
a & b & c & d \\ \hline
7 & y & 4 & 17 \\
7 & y & 6 & 20 \\
7 & y & 8 & 22 \\
7 & y & 10 & 23 \\
9 & z & 4 & 18 \\
9 & z & 6 & 21 \\
9 & z & 8 & 23 \\
9 & z & 10 & 24 
\end{array}
\qqqquad
\vcenter{\xymatrix@M=1ex@C=20ex@R=10ex@!0{
& \alpha(\gamma(f)) ;_\gamma \beta(g)
\ar@{=>}[dl]_{(a)} \ar@{|=|}[dr]^{(b);1} \\
\alpha(\gamma(f) ;_\gamma \beta(g)) \ar@{=>}[ddd]_{\alpha((c))}
&& \gamma(\alpha(f)) ;_\gamma \beta(g) \ar@{=>}[d]^{(c)} \\
&& \beta(\gamma(\alpha(f)) ;_\gamma g) \ar@{|=|}[d]^{\beta((b);1)} \\
&& \beta(\alpha(\gamma(f)) ;_\gamma g) \ar@{=>}[d]^{\beta((a))} \\
\alpha(\beta(\gamma(f) ;_\gamma g)) \ar@{|=|}[rr]_{(d)}
&& \beta(\alpha(\gamma(f) ;_\gamma g))}}
\]
where $w \in \{15, 16, 17, 18\}$, $x \in \{16, 19, 20, 21\}$, $y \in
\{17,20,22,23\}$, and $z \in \{18,21,23,24\}$. This means that there are 64
reductions that fit this general case!

\item[Reduction diagram 5:] each row in the table corresponds to the
resolution of the critical pair ($a$)-($b$).
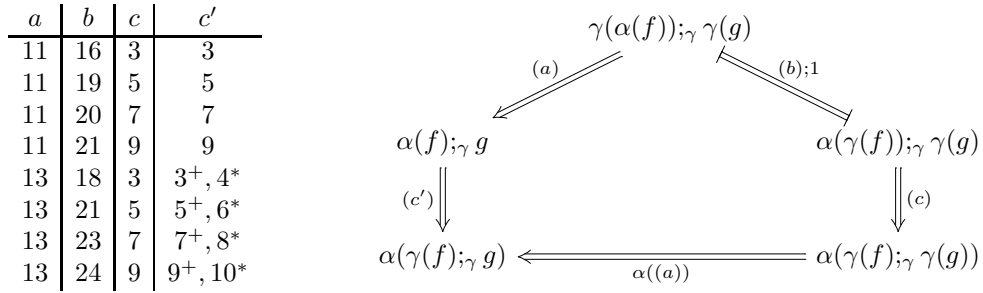
\[\begin{array}{c|c|c|c}
a & b & c & c'\\ \hline
11 & 16 & 3 & 3 \\
11 & 19 & 5 & 5 \\
11 & 20 & 7 & 7 \\
11 & 21 & 9 & 9 \\
13 & 18 & 3 & 3^+,4^* \\
13 & 21 & 5 & 5^+,6^* \\
13 & 23 & 7 & 7^+,8^* \\
13 & 24 & 9 & 9^+,10^*
\end{array}
\qqqquad
\vcenter{\xymatrix@M=1ex@C=20ex@R=10ex@!0{
& \gamma(\alpha(f)) ;_\gamma \gamma(g)
\ar@{=>}[dl]_{(a)} \ar@{|=|}[dr]^{(b);1} \\
\alpha(f) ;_\gamma g \ar@{=>}[d]_{(c')}
&& \alpha(\gamma(f)) ;_\gamma \gamma(g) \ar@{=>}[d]^{(c)} \\
\alpha(\gamma(f) ;_\gamma g) && \alpha(\gamma(f) ;_\gamma \gamma(g))
\ar@{=>}[ll]^-{\alpha((a))}}}
\]
where $x^+,y^*$ means zero or one application of the rewrite ($x$) and zero
or more applications of the rewrite ($y$).

This may be a good time for a concrete example. Suppose $(a,b,c,c') =
(13,21,5,(5^+,6^*))$, $i \in \{1,\ldots,n\}$, and $\alpha \in \Omega_k$.
The reduction diagram for this case is: 
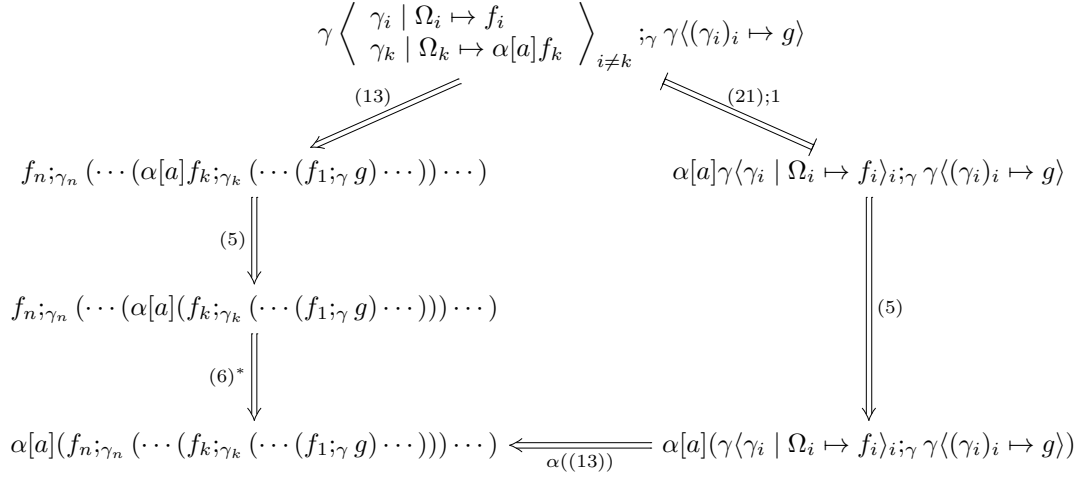
\[\xymatrix@M=1ex@C=27ex@R=12ex@!0{
& {\gamma\left\<\begin{array}{l}
\gamma_i \mid \Omega_i \mapsto f_i \\
\gamma_k \mid \Omega_k \mapsto \alpha[a]f_k
\end{array}\right\>_{i \neq k} ;_\gamma \gamma\<(\gamma_i)_i \mapsto g\>}
\ar@{=>}[dl]_-{(13)} \ar@{|=|}[dr]^-{(21);1} \\
f_n ;_{\gamma_n} (\cdots (\alpha[a]f_k ;_{\gamma_k} (\cdots (f_1
;_\gamma g)\cdots))\cdots) \ar@{=>}[d]_{(5)}
&& \alpha[a]\gamma\<\gamma_i \mid \Omega_i \mapsto f_i\>_i ;_\gamma
\gamma\<(\gamma_i)_i \mapsto g\> \ar@{=>}[dd]^{(5)} \\
f_n ;_{\gamma_n} (\cdots (\alpha[a](f_k ;_{\gamma_k} (\cdots (f_1
;_\gamma g)\cdots)))\cdots) \ar@{=>}[d]_{(6)^*} \\
\alpha[a](f_n ;_{\gamma_n} (\cdots (f_k ;_{\gamma_k} (\cdots (f_1
;_\gamma g)\cdots)))\cdots)
&& \alpha[a](\gamma\<\gamma_i \mid \Omega_i \mapsto f_i\>_i ;_\gamma
\gamma\<(\gamma_i)_i \mapsto g\>)
\ar@{=>}[ll]^-{\alpha((13))}}
\]

\item[Reduction diagram 6:] each row in the table corresponds to the
resolution of the critical pair ($a$)-($b$).
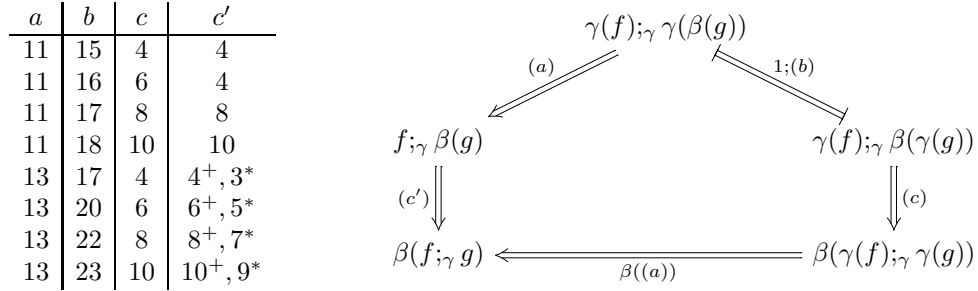
\[\begin{array}{c|c|c|c}
a & b & c & c' \\
\hline
11 & 15 & 4 & 4 \\
11 & 16 & 6 & 4 \\
11 & 17 & 8 & 8 \\
11 & 18 & 10 & 10 \\
13 & 17 & 4 & 4^+,3^* \\
13 & 20 & 6 & 6^+,5^* \\
13 & 22 & 8 & 8^+,7^* \\
13 & 23 & 10 & 10^+,9^*
\end{array}
\qqqquad
\vcenter{\xymatrix@M=1ex@C=20ex@R=10ex@!0{
& \gamma(f) ;_\gamma \gamma(\beta(g))
\ar@{=>}[dl]_{(a)} \ar@{|=|}[dr]^{1;(b)} \\
f ;_\gamma \beta(g) \ar@{=>}[d]_{(c')}
&& \gamma(f) ;_\gamma \beta(\gamma(g)) \ar@{=>}[d]^{(c)} \\
\beta(f ;_\gamma g) && \beta(\gamma(f) ;_\gamma \gamma(g))
\ar@{=>}[ll]^-{\beta((a))}}}
\]
where $x^+,y^*$ means zero or one application of the rewrite ($x$) and zero
or more applications of the rewrite ($y$).
\end{description}

We now explore the cases when the index sets may be empty. For the empty
ltensor (dually rpar) rule the rewrites fit the cases above. The cases for
the empty cotuple (dually tuple) and rtensor (dually lpar) however do not.
We start by first examing what happens to the reduction diagrams for the
case of the empty cotuple (dually tuple).

\begin{description}
\item[Reduction diagram 1:]
$\xymatrix{\alpha\{\,\};1 \ar@{=>}[r]_{(3)}^{(1)} & \alpha\{\,\}}$

\item[Reduction diagram 2:] there are two (non-dual) cases corresponding
to only $I = \emptyset$ and both $I = J = \emptyset$. Each row in the table
corresponds to the resolution of the critical pair (3)-($a$).
\[\begin{array}{c|c} 
a & b \\ \hline
4 & 15 \\
6 & 16 \\
8 & 17 \\
10 & 18
\end{array}
\qqqquad
\vcenter{\xymatrix@M=1ex@C=12ex@R=9ex@!0{
& \alpha\{\,\} ;_\gamma \beta(g) \ar@{=>}[dl]_{(3)} \ar@{=>}[dr]^{(a)} \\
\alpha\{\,\} \ar@{|=|}[dr]_-{(b)}
&& \beta(\alpha\{\,\} ;_\gamma g) \ar@{=>}[dl]^{\beta((3))} \\
& \beta(\alpha\{\,\})}}
\]
Dual to the above diagram is the nullary redutions for $\beta$. If both
$\alpha$ and $\beta$ have empty index sets:
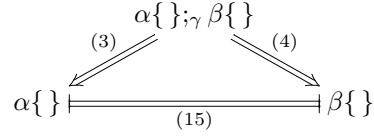
\[\xymatrix@R=4ex@C=5ex{
& \alpha\{\,\} ;_\gamma \beta\{\,\} \ar@{=>}[dl]_{(3)} \ar@{=>}[dr]^{(4)} \\
\alpha\{\,\} \ar@{|=|}[rr]_-{(15)} && \beta\{\,\}}
\]

\item[Reduction diagram 3:] there are three cases. The first we describe is
when both $I = J = \emptyset$. The resolution is as follows:
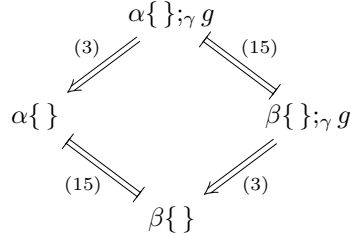
\[\xymatrix@M=1ex@C=12ex@R=9ex@!0{
& \alpha\{\,\} ;_\gamma g \ar@{=>}[dl]_{(3)} \ar@{|=|}[dr]^{(15)} \\
\alpha\{\,\} \ar@{|=|}[dr]_-{(15)}
&& \beta\{\,\} ;_\gamma g \ar@{=>}[dl]^{(3)} \\
& \beta\{\,\}}
\]
The two remaining cases corresponding to whether the apex (of the reduction
diagram) starts with $\alpha\{\,\}$ or with $\beta(\alpha\{\,\})$. Each row
in the table corresponds to the resolution of the critical pair (3)-($b$) in
the reduction diagram on the left and ($a$)-($b$) in the reduction diagram on
the right.
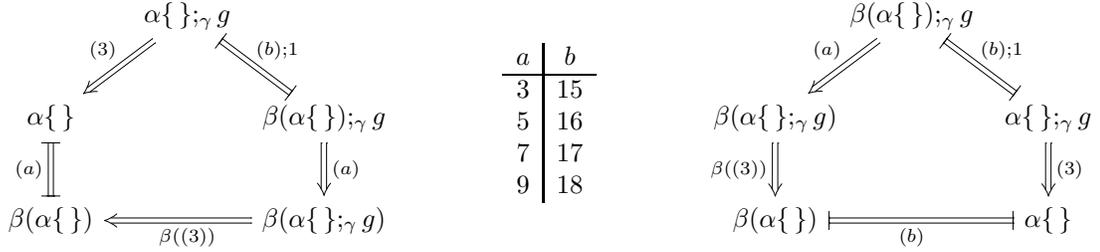
\[\vcenter{\xymatrix@M=1ex@C=12ex@R=9ex@!0{
& \alpha\{\,\} ;_\gamma g \ar@{=>}[dl]_{(3)} \ar@{|=|}[dr]^{(b);1} \\
\alpha\{\,\} \ar@{|=|}[d]_{(a)}
&& \beta(\alpha\{\,\}) ;_\gamma g \ar@{=>}[d]^{(a)} \\
\beta(\alpha\{\,\})
&& \beta(\alpha\{\,\} ;_\gamma g) \ar@{=>}[ll]^{\beta((3))}}}
\qqqquad
\begin{array}{c|c} 
a & b \\ \hline
3 & 15 \\
5 & 16 \\
7 & 17 \\
9 & 18
\end{array}
\qqqquad
\vcenter{\xymatrix@M=1ex@C=12ex@R=9ex@!0{
& \beta(\alpha\{\,\}) ;_\gamma g \ar@{=>}[dl]_{(a)}
\ar@{|=|}[dr]^{(b);1} \\
\beta(\alpha\{\,\} ;_\gamma g) \ar@{=>}[d]_{\beta((3))}
&& \alpha\{\,\} ;_\gamma g \ar@{=>}[d]^{(3)} \\
\beta(\alpha\{\,\}) && \alpha\{\,\} \ar@{|=|}[ll]^{(b)}}}
\]

\item[Reduction diagram 4:] The channel $\alpha$ is non-empty and if channel
$\gamma$ is empty the reduction diagram is identical.

\item[Reduction diagram 5:] each row in the table corresponds to the
resolution of the critical pair ($a$)-($b$).
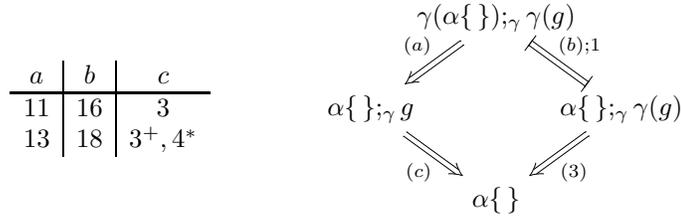
\[\begin{array}{c|c|c}
a & b & c \\ \hline
11 & 16 & 3 \\
13 & 18 & 3^+,4^*
\end{array}
\qqqquad
\vcenter{\xymatrix@M=1ex@C=11ex@R=8ex@!0{
& \gamma(\alpha\{\,\}) ;_\gamma \gamma(g)
\ar@{=>}[dl]_-{(a)} \ar@{|=|}[dr]^-{(b);1} \\
\alpha\{\,\} ;_\gamma g \ar@{=>}[dr]_-{(c)}
&& \alpha\{\,\} ;_\gamma \gamma(g) \ar@{=>}[dl]^-{(3)} \\
& \alpha\{\,\}}}
\]
where $3^+,4^*$ means zero or one application of the rewrite (3) and zero
or more applications of the rewrite (4).

\item[Reduction diagram 6:] each row in the table corresponds to the
resolution of the critical pair ($a$)-($b$).
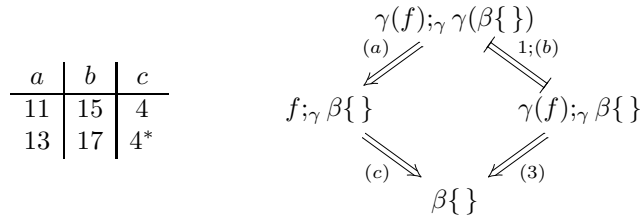
\[\begin{array}{c|c|c}
a & b & c \\ \hline
11 & 15 & 4 \\
13 & 17 & 4^*
\end{array}
\qqqquad
\vcenter{\xymatrix@M=1ex@C=11ex@R=8ex@!0{
& \gamma(f) ;_\gamma \gamma(\beta\{\,\}) \ar@{=>}[dl]_-{(a)}
 \ar@{|=|}[dr]^-{1;(b)} \\
f ;_\gamma \beta\{\,\} \ar@{=>}[dr]_-{(c)}
&& \gamma(f) ;_\gamma \beta\{\,\} \ar@{=>}[dl]^-{(3)} \\
& \beta\{\,\}}}
\]
where $4^*$ means zero or more applications of the rewrite (4).

\end{description}

We now examing what happens to the reduction diagrams for the case of the
empty rtensor (dually lpar). Due to the typing constraints the only
reduction diagrams that need be considered are 1 and 6.

\begin{description}
\item[Reduction diagram 1:]
$\xymatrix{\alpha\<\,\>;1 \ar@{=>}[r]_{(13)}^{(1)} & \alpha\<\,\>}$

\item[Reduction diagram 6:] each row in the table corresponds to the
resolution of the critical pair (13)-($a$).
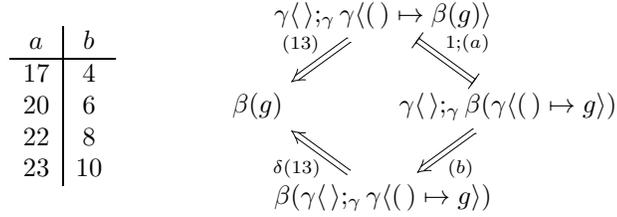
\[\begin{array}{c|c}
a & b \\ \hline
17 & 4 \\
20 & 6 \\
22 & 8 \\
23 & 10 
\end{array}
\qqqquad
\vcenter{\xymatrix@M=1ex@C=11ex@R=8ex@!0{
& \gamma\<\,\> ;_\gamma \gamma\<(\,) \mapsto \beta(g)\>
\ar@{=>}[dl]_-{(13)} \ar@{|=|}[dr]^-{1;(a)} \\
\beta(g)
&& \gamma\<\,\> ;_\gamma \beta(\gamma\<(\,) \mapsto g\>) \ar@{=>}[dl]^-{(b)} \\
& \beta(\gamma\<\,\> ;_\gamma \gamma\<(\,) \mapsto g\>)
\ar@{=>}[ul]^-{\delta(13)}}}
\]
\end{description}

\subsubsection{The measure on the rewriting arrows}

We define a measure $\lambda:A \ra \bag(\mathbb{N})$ on the rewriting arrows
as follows:
\begin{itemize}
\item if $\xymatrix{t_1 \ar@{=>}[r]^x & t_2}$ then $\lambda(x) =
\min\{\Lambda(t_1),\ \Lambda(t_2)\}$
\item if $\xymatrix{t_1 \ar@{|=|}[r]^x & t_2}$ then $\lambda(x) =
\max\{\Lambda(t_1),\ \Lambda(t_2)\}$
\end{itemize}
where $\Lambda(t)$ is the bag of cut heights of $t$.

A quick examination of the reduction diagrams now confirms that this measure
will decrease when we replace a divergence with a convergence.

This completes the proof of the proposition:

\begin{proposition}
$\channel(\cat{A})$ under the rewrites (1)-(14) is confluent modulo the
equations (15)-(24).
\end{proposition}

\section{Proof that $\channel(\cat{A})$ is a polycategory}
\label{sec-proof-polycat}

In order to show that $\channel(\cat{A})$ is a polycategory three
properties must be satisified:
\begin{enumerate}[{\upshape (i)}]
\item $\channel(\cat{A})$ must have identity maps which behave in the
correct manner,
\item composition in $\channel(\cat{A})$ must be associative, and
\item composition in $\channel(\cat{A})$ must satisfy the interchange law.
\end{enumerate}

We prove each in turn.

\begin{lemma}
The identity acts as a neutral element with respect to composition. That is,
given terms of the form
\[f:: \Gamma \ra \Delta,\gamma \cc X \quad \text{and} \quad 1_X::\gamma \cc
X \ra \delta \cc X
\]
we have $f ;_\gamma 1_X = f$ (upto renaming channels), and dually, given terms
of the form
\[1_X::\gamma \cc X \ra \delta \cc X \quad \text{and} \quad f:: \delta \cc
X,\Gamma \ra \Delta
\]
we have $1_X ;_\delta f =f$ (upto renaming channels).
\end{lemma}

\begin{proof}
We shall suppose the identity is on the left; duality covers the other case.
The case where $f$ is the identity is clearly true. So suppose $f$ is of the
form $\alpha(f)$, where $\alpha \neq \gamma$. If $\alpha(f)$ is a tuple or
rtensor with $I = \emptyset$ then $\alpha(f) ;_\gamma 1_X$ will reduce to
$\alpha\{\,\}$ or $\delta\<\,\>$ respectively. Otherwise, any rewrite
applied to $\alpha(f) ;_\gamma 1_X$ will move the identity in $\alpha(f
;_\gamma 1_X)$ moving the cut onto a smaller term from which we may apply
the inductive hypothesis. Thus, the only cases we must explore is when the
term $f$ operates on $\gamma$ at the top level.

The proof will proceed by structure induction on the term $\gamma(f)$.
Without loss of generality we may assume that $f$ is cut free.

\begin{itemize}
\item The base case is a cut with an atomic identity: here the
cut-elimination step removes the identity and the result is immediate.

\item $f = \gamma\{a_i \mapsto f_i\}_{i \in I}$. There are two cases to
consider corresponding to $I = \emptyset$ and $I \neq \emptyset$.
In the first case $1_X = \delta\{\,\}$ and by rewrite (4) this reduces
to $\delta\{\,\}$. In the second case we have $X=\prod_i a_i:X_i$ and
$1_X = \delta\{a_i \mapsto \gamma[a_i] 1_{X_i}\}_i$. This gives
\[\begin{array}{rcl}
\gamma\{a_i \mapsto f_i\}_i ;_\gamma \delta\{a_i \mapsto \gamma[a_i]
1_{X_i}\}_i
&\Lra& \delta\{a_i \mapsto \gamma\{a_i \mapsto f_i\}_i ;_\gamma
\gamma[a_i] 1_{X_i}\}_i \medskip\\
&\Lra& \delta\{a_i \mapsto f_i ;_\gamma 1_{X_i}\}_i
\end{array}
\]
which moves the composition onto smaller terms. Applying the inductive
hypothesis now gives the desired result.

\item $f = \gamma[a_k]f'$. In this case $X= \sum_i X_i$ and $1_X =
\gamma\{a_i \mapsto \delta[a_i] 1_{X_i}\}_i$ which gives
\[\begin{array}{rcl}
\gamma[a_k]f' ;_\gamma \gamma\{a_i \mapsto \delta[a_i] 1_{X_i}\}_i
&\Lra& f' ;_\gamma \delta[a_k] 1_{X_k} \medskip\\
&\Lra& \delta[a_k](f' ;_\gamma 1_{X_k})
\end{array}
\]
which moves the composition onto a smaller term. Applying the inductive
hypothesis now gives the desired result.

\item $f = \gamma\<(\gamma_i)_i \mapsto f'\>$. In this case $X = \bigot_i
\gamma_i:X_i$ (in $f$ and the domain of $1_X$ and $X = \bigot_i \delta_i:X_i$
in the codomain of $1_X$) and $1_X = \delta\<(\delta_i)_i \mapsto
\gamma\<\gamma_i \mid \delta_i \mapsto 1_{X_i}\>_i\>$. In this case suppose
$I = \{1,\ldots,n\}$. This gives
\[\begin{array}{rcl}
\gamma\<(\gamma_i)_i \mapsto f'\> ;_\gamma \delta\<(\delta_i)_i \mapsto
\gamma\<\gamma_i \mid \delta_i \mapsto 1_{X_i}\>_i\>
&\Lra&
\delta\<(\delta_i)_i \mapsto \gamma\<(\gamma_i)_i \mapsto f'\> ;_\gamma 
\gamma\<\gamma_i \mid \delta_i \mapsto 1_{X_i}\>_i\> \medskip\\
&\Lra&
\delta\<(\delta_i)_i \mapsto ( \cdots (f' ;_{\gamma_n} 1_{X_n}) \cdots )
;_{\gamma_1} 1_{X_1}\>
\end{array}
\]
which moves the composition onto smaller terms. Applying the inductive
hypothesis now gives the desired result.

\item $f = \gamma\<\gamma_i \mid \Omega_i \mapsto f_i\>_{i \in I}$. There
are two cases to consider corresponding to $I = \emptyset$ and $I \neq
\emptyset$. In the first case $1_X = \gamma\<(\,) \mapsto \delta\<\,\>$ which
by (13) reduces to $\delta\<\,\>$. In the second case we have $X = \bigox_i
\gamma_i:X_i$ (in $f$ and the domain of $1_X$ and $X = \bigox_i \delta_i:X_i$
in the codomain of $1_X$) and $1_X = \gamma\<(\gamma_i)_i \mapsto
\delta\<\delta_i \mid \gamma_i \mapsto 1_{X_i}\>_i\>$. In this case suppose
$I = \{1,\ldots,n\}$. This gives
\[\begin{array}{rcl}
\gamma\<\gamma_i \mid \Omega_i \mapsto f_i\>_i ;_\gamma \gamma\<
(\gamma_i)_i \mapsto \delta\<\delta_i \mid \gamma_i \mapsto 1_{X_i}\>_i\>
&\Lra&
f_n ;_{\gamma_n}(\cdots (f_1 ;_{\gamma_1} \delta\<\delta_i \mid \gamma_i
\mapsto 1_{X_i}\>_i) \cdots ) \medskip\\
&\Lra& \delta\<\delta_i \mid \Omega_i \mapsto f_i ;_{\gamma_i} 1_{X_i}\>_i
\end{array}
\]
which moves the composition onto smaller terms. Applying the inductive
hypothesis now gives the desired result.
\end{itemize}

This now completes the proof that the identity acts as a neutral element with
respect to composition in this system.
\end{proof}

\begin{lemma}
Cut satisfies the associative law. That is, given terms of the form
\[f::\Gamma_1 \ra \Gamma_2,\gamma \cc X \qquad 
g::\gamma \cc X,\Delta_1 \ra \Delta_2,\delta \cc Y \qquad 
h::\delta \cc Y,\Phi_1 \ra \Phi_2
\]
the composities $(f ;_\gamma g) ;_\delta h$ and $f ;_\gamma (g ;_\delta h)$
are $\sim$-equivalent.
\end{lemma}

\begin{proof}
The proof proceeds by structural induction on $f$, $g$, and $h$; without
loss of generality we may assume that $f$, $g$, and $h$ are cut free.
Without explicitly mentioning where, duality will be used to reduce the
number of cases presented. To show that the composites are $\sim$-equivalent
we will make use of the generalized rewrites. Additionally, the notation
\[\xymatrix{f ;_\gamma g ;_\delta h \ar@{=>}[r]^{;_l} & (f ;_\gamma g)
;_\delta h} \quad \text{and} \quad \xymatrix{
f ;_\gamma g ;_\delta h \ar@{=>}[r]^{;_r} & f ;_\gamma (g ;_\delta h)}
\]
will be used to indicate the two possible composites of
$f ;_\gamma g ;_\delta h$. 

Recall that if two terms are $\sim$-equavalent then they must be related
through the communting conversions alone.

If $f = 1_X$ then $\xymatrix{1_X ;_\gamma g ;_\delta h
\ar@{=>}[r]^-{;_l}_-{;_r} & g ;_\delta h}$ shows that the composites are
$\sim$-equivalent. 

Suppose now that $f$ is of the form $\alpha(f)$, where $\alpha \neq \gamma$.
If $\alpha(f)$ is a tuple or rtensor with $I = \emptyset$ then $f ;_\gamma 
g ;_\delta h$ will respectively reduce to $\alpha\{\,\}$ or $\beta\<\,\>$
(where $\beta$ is the codomain channel of $h$). If $I \neq \emptyset$ then
in following reduction diagram
\[\xymatrix@M=1ex@C=12ex@R=9ex@!0{
& \alpha(f) ;_\gamma g ;_\delta h \ar@{=>}[dl]_{;_l} \ar@{=>}[dr]^{;_r} \\
\alpha(f ;_\gamma g) ;_\delta h \ar@{=>}[d]_{;} 
&& \alpha(f) ;_\gamma (g ;_\delta h) \ar@{=>}[d]^{;}  \\
\alpha((f ;_\gamma g) ;_\delta h) \ar@{|=|}[rr]_{*}
&& \alpha(f ;_\gamma (g ;_\delta h))}
\]
the composites are moved onto smaller terms, which by induction are
satisfy the associative law, and hence, both composites are
$\sim$-equavalent.

It remains to examine the case where $f$ is of the form $\gamma(f)$. There
are two cases corresponding to whether $g$ is of the form $\beta(g)$, where
$\beta \neq \gamma$ and $\beta \neq \delta$, or $\gamma(g)$. (The dual
rewrites will suffice for the $\beta = \delta$ case.) In the first
case if $\beta(g) = \beta\{\,\}$ then both composites will reduce to
$\beta\{\,\}$. Note the $\beta(g) \neq \beta\<\,\>$ as both composites are
defined. So suppose that $\beta(g) \neq \beta\{\,\}$. The following reduction
diagram shows that both composites are $\sim$-equivalent (by induction):
\[\xymatrix@M=1ex@C=12ex@R=9ex@!0{
& \gamma(f) ;_\gamma \beta(g) ;_\delta h \ar@{=>}[dl]_{;_l}
\ar@{=>}[dr]^{;_r} \\
\beta(\gamma(f) ;_\gamma g) ;_\delta h \ar@{=>}[d]_{;} 
&& \gamma(f) ;_\gamma \beta(g ;_\delta h) \ar@{=>}[d]^{;}  \\
\beta((\gamma(f) ;_\gamma g) ;_\delta h) \ar@{|=|}[rr]_{*}
&& \beta(\gamma(f) ;_\gamma (g ;_\delta h))}
\]

In the second case $(\beta = \gamma)$ there is no need to generalize as
there are only two (non-dual) rewrites (which also has the benefit of
providing a ``concrete'' example). The first is when the left-hand side
composite is the rewrite (11):
\[\xymatrix@M=1ex@C=25ex@R=9ex@!0{
\gamma\{a_i \mapsto f_i\}_i ;_\gamma \gamma[a_k]g ;_\delta h
\ar@{=>}[r]^-{;_l} \ar@{=>}[d]_-{;_r} &
(f_k ;_\gamma g) ;_\delta h \\
\gamma\{a_i \mapsto f_i\}_i ;_\gamma \gamma[a_k](g ;_\delta h) 
\ar@{=>}[r]_-{;} & f_k ;_\gamma (g ;_\delta h) \ar@{|=|}[u]_{*}}
\]

The second case is when the left-hand side composite is the rewrite (13). In
this case assume that $I = \{1,\ldots,n\}$:
\[\xymatrix@M=1ex@C=40ex@R=12ex@!0{
\alpha\<\alpha_i \mid \Omega_i \mapsto f_i\>_i ;_\gamma \gamma\<(\gamma_i)
\mapsto g\> ;_\delta h
\ar@{=>}[d]_-{;_r} \ar@{=>}[r]^-{;_l}
& (f_n ;_{\gamma_n} (\cdots (f_1 ;_{\gamma_1} g) \cdots)) ;_\delta h
\ar@{|=|}[d]^{*} \\
\alpha\<\alpha_i \mid \Omega_i \mapsto f_i\>_i ;_\gamma \gamma\<(\gamma_i)
\mapsto g ;_\delta h\>
\ar@{=>}[r]_-{;} & 
(f_n ;_{\gamma_n} (\cdots (f_1 ;_{\gamma_1} (g ;_\delta h)) \cdots))}
\]

If all of $f$, $g$, and $h$ are atomic then composition is associative as it
must be associative in the underlying polycategory. If some of $f$, $g$,
and $h$ are atomic a quick check of the possibilities will show that one
ends up with a case similar to one of the cases above.

Thus, composition is associative.
\end{proof}

\begin{lemma}
Cut satisfies the interchange property. That is, given terms of the form
\[f:: \Gamma_1 \ra \gamma \cc X,\Gamma_2,\delta \cc Y \qquad g::
\Delta_1,\gamma \cc X \ra \Delta_2 \qquad h::\delta \cc Y,\Phi_1 \ra \Phi_2
\]
the composites $(f ;_\gamma g) ;_\delta h$ and $(f ;_\delta h) ;_\gamma g$
are $\sim$-equivalent. Dually, given sequents of the form
\[f:: \Gamma_1 \ra \Gamma_2,\gamma \cc X \qquad g::\Delta_1 \ra \delta \cc
Y,\Delta_2 \qquad h::\gamma \cc X,\Phi_1,\delta \cc Y \ra \Phi_2
\]
the composites $f ;_\gamma (g ;_\delta h)$ and $g ;_\delta (f ;_\gamma h)$
are $\sim$-equivalent.
\end{lemma}

\begin{proof}
The proof proceeds by structural induction on $f$, $g$, and $h$; without
loss of generality we may assume that $f$, $g$, and $h$ are cut free.
Without explicitly mentioning where, duality will be used to reduce the
number of cases presented. In the following
\[\xymatrix{f ;_\gamma g ;_\delta h \ar@{=>}[r]^{;_g} & (f ;_\gamma g)
;_\delta h} \quad \text{and} \quad \xymatrix{f ;_\gamma g ;_\delta h
\ar@{=>}[r]^{;_h} & (f ;_\delta h) ;_\gamma g}
\]
will be used respectively to indicate composing first with $g$ and composing
first with $h$.

Notice that $f$ may not be $1_X$ as it requires at least two codomain
channels. So suppose that $f$ is of the form $\alpha(f)$. If $f =
\alpha\{\,\}$ then both composites will reduce to $\alpha\{\,\}$. The term
$f$ may not be the empty rtensor since, as noted above, it requires at least
two codomain channels. If $I \neq \emptyset$ then in the reduction diagram
\[\xymatrix@M=1ex@C=12ex@R=9ex@!0{
& \alpha(f) ;_\gamma g ;_\delta h \ar@{=>}[dl]_{;_g} \ar@{=>}[dr]^{;_h} \\
\alpha(f ;_\gamma g) ;_\delta h \ar@{=>}[d]_{;} 
&& \alpha(f ;_\delta h) ;_\gamma g \ar@{=>}[d]^{;}  \\
\alpha((f ;_\gamma g) ;_\delta h) \ar@{|=|}[rr]_{*}
&& \alpha((f ;_\delta h) ;_\gamma g)}
\]
the composites are moved onto smaller terms, which by induction satisfy the
interchange law, and hence, both composites are $\sim$-equivalent.

It remains to examine the case where $f$ is of the form $\gamma(f)$. There
are two cases corresponding to whether $g$ is of the form $\beta(g)$, where
$\beta \neq \gamma$ and $\beta \neq \delta$, or $\gamma(g)$. (The dual
rewrites will suffice for the $\beta = \delta$ case.) In the first case if
$\beta(g) = \beta\{\,\}$ then both composites will reduce to $\beta\{\,\}$.
If $\beta(g) = \beta\<\,\>$ (empty lpar) then both composites will reduce to
$f ;_\delta h$.

Now suppose that $\beta(g)$ is not a nullary operation. The following
reduction diagram shows that both composites are $\sim$-equivalent (by
induction):
\[\xymatrix@M=1ex@C=12ex@R=9ex@!0{
& \gamma(f) ;_\gamma \beta(g) ;_\delta h \ar@{=>}[dl]_{;_g}
\ar@{=>}[dr]^{;_h} \\
\beta(\gamma(f) ;_\gamma g) ;_\delta h \ar@{=>}[d]_{;} 
&& \gamma(f ;_\delta h) ;_\gamma \beta(g) \ar@{=>}[d]^{;}  \\
\beta((\gamma(f) ;_\gamma g) ;_\delta h) \ar@{|=|}[rr]_{*}
&& \beta(\gamma(f ;_\delta h) ;_\gamma g)}
\]

In the second case $(\beta=\gamma)$ there is only one possible choice:
composition with $g$ is an application of the rewrite (11). Explicitly,
\[\xymatrix@M=1ex@C=25ex@R=9ex@!0{
\gamma\{a_i \mapsto f_i\}_i ;_\gamma \gamma[a_k]g ;_\delta h
\ar@{=>}[r]^-{;_l} \ar@{=>}[d]_-{;_h} &
(f_k ;_\gamma g) ;_\delta h \\
\gamma\{a_i \mapsto f_i ;_\delta h\}_i ;_\gamma \gamma[a_k]g 
\ar@{=>}[r]_-{;} & (f_k ;_\delta h) ;_\gamma g \ar@{|=|}[u]_{*}}
\]
which shows that the composites are $\sim$-equivlant.

If all of $f$, $g$, and $h$ are atomic then composition satisfies the
interchange law as it must satisfy the interchange law in the underlying
polycategory. If some of $f$, $g$, and $h$ are atomic a quick check of the
possibilities will show that one ends up with a case similar to one of the
cases above.

Thus, composition satifies the interchange law.
\end{proof}

\end{document}